\documentclass[11pt,a4paper,leqno,numbers=noenddot,abstract=false]{scrartcl}

\setkomafont{disposition}{\normalfont\bfseries} 

\usepackage[english]{babel}
\usepackage[T1]{fontenc}
\usepackage[utf8]{inputenc}

\usepackage{amssymb}
\usepackage{amsmath}
\usepackage{lmodern}
\usepackage{old-arrows}

\usepackage[colorlinks=true,allcolors=blue,bookmarks=true,linktocpage=true]{hyperref}

\makeatletter
\newcommand{\labeltext}[2]{%
  \@bsphack
  \csname phantomsection\endcsname 
  \def\@currentlabel{#1}{\label{#2}}%
  \@esphack
} \makeatother 

\usepackage{enumitem}
\usepackage{mathtools} 
\usepackage{graphicx}
\usepackage[format=plain]{caption} 
\usepackage[labelformat=simple]{subcaption}

\usepackage[style=authoryear-comp,maxnames=5,maxcitenames=2,dashed=false,giveninits=true,
            useprefix=true]{biblatex}
\DeclareFieldFormat 
  [article,inbook,incollection,inproceedings,patent,thesis,unpublished]
  {title}{#1\isdot} 
\DeclareNameAlias{sortname}{family-given} 
\DeclareFieldFormat{pages}{#1} 
\renewbibmacro{in:}{%
  \ifentrytype{article}{}{\printtext{\bibstring{in}\intitlepunct}}}
\renewbibmacro*{volume+number+eid}{
  \printfield{volume}%
  \printfield{number}%
  \setunit{\addcomma\space}%
  \printfield{eid}}
\DeclareFieldFormat[article,periodical]{volume}{\mkbibbold{#1}} 
\DeclareFieldFormat{doi}{%
  \hfil\penalty90\hfilneg\space \textsc{doi}\addcolon\addnbspace
  \ifhyperref
    {\href{https://doi.org/#1}{\nolinkurl{#1}}}
    {\nolinkurl{#1}}}
\DefineBibliographyExtras{english}{%
  
} 

\renewcommand{\cite}{\textcite}

\newcommand{\citelink}[2]{\hyperlink{cite.\therefsection @#1}{#2}}
\newcommand{\HMF}{\citelink{Olver2010}{HMF}}

\DefineBibliographyStrings{english}{
  andothers = {\emph{et al.}},
} 

 \DeclarePunctuationPairs{colon}{*}

\DeclareSortingNamekeyTemplate{ 
  \keypart{
    \namepart{family}
  }
  \keypart{
    \namepart{prefix}
  }
  \keypart{
    \namepart{given}
  }
  \keypart{
    \namepart{suffix}
  }
} 

\renewbibmacro*{begentry}{\midsentence}
\AtBeginDocument{\toggletrue{blx@useprefix}} 

\usepackage[mathscr]{euscript}
\let\mathcal\mathscr

\newcommand{\Z}{\mathbb{Z}}
\newcommand{\R}{\mathbb{R}}
\newcommand{\C}{\mathbb{C}}
\newcommand{\dss}{\displaystyle}
\newcommand{\bfP}{\mathbf{P}}
\newcommand{\bft}{\mathbf{t}}
\newcommand{\bfd}{\mathbf{d}}
\renewcommand{\theta}{\vartheta}
\renewcommand{\epsilon}{\varepsilon}

\DeclarePairedDelimiter\abs{\lvert}{\rvert} 

\DeclareMathOperator{\RE}{Re} \DeclareMathOperator{\Log}{Log} 

\DeclareMathAlphabet{\mathmybb}{U}{bbold}{m}{n}
\newcommand{\1}{\mathmybb{1}}

\usepackage[amsmath,amsthm,thmmarks]{ntheorem}

\makeatletter 
 \newtheoremstyle{mychange}%
  {\item[\hskip\labelsep \theorem@headerfont (##2)\ ##1\theorem@separator]}%
  {\item[\hskip\labelsep \theorem@headerfont (##2)\ ##1\ (##3)\theorem@separator]}
\makeatother

\theoremstyle{mychange}%
\theorembodyfont{\itshape}%
\theoremheaderfont{\upshape}%
\theoremnumbering{Alph}%
\newtheorem{regtheorem}{}

\theoremstyle{remark}
\newtheorem*{remark}{Remark}

\newcommand{\qedhere}{\null}

\renewcommand*{\thefootnote}{(\arabic{footnote})}

\let\svthefootnote\thefootnote
\newcommand\freefootnote[1]{%
  \let\thefootnote\relax%
  \footnotetext{#1}%
  \let\thefootnote\svthefootnote%
}

\allowdisplaybreaks[4] 

\begin{document}

\title{On the number of elements beyond the ones actually observed}
\author{by Eugenio Regazzini}
\date{\normalsize\itshape Dipartimento di Matematica, Università degli Studi di Pavia\\ 
\href{mailto:eugenio.regazzini@unipv.it}{eugenio.regazzini@unipv.it}} 

\maketitle

\begin{abstract}\small
  In this work, a variant of the birth and death chain with constant intensities, originally introduced by 
  Bruno de Finetti way back in 1957, is revisited. This fact is also underlined by the choice of the title, 
  which is clearly a literal translation of the original one. Characteristic of the variant is that it 
  allows negative jumps of any magnitude. And this, as explained in the paper, might be useful in offering 
  some insight into the issue, arising in numerous situations, of inferring the number of the undetected 
  elements of a given population. One thinks, for example, of problems concerning abundance or richness of 
  species.
  
  The author's purpose is twofold: to align the original de Finetti's construction with the modern, 
  well-established theory of the continuous-time Markov chains with discrete state space and show how it 
  could be used to make probabilistic previsions  on the number of the unseen elements of a population. 
  With the aim of enhancing the possible practical applications of the model, one discusses the statistical 
  point estimation of the rates which characterize its infinitesimal description. 
\end{abstract}

\freefootnote{\href{https://mathscinet.ams.org/mathscinet/msc/msc2020.html}{\textit{MSC2020 subject 
classifications.}} Primary 60J27, 60J28, secondary 62F10.} \freefootnote{\textit{Key words and phrases.} 
Birth-death Markov chain, continuous-time Markov process, estimation by the method of moments, incomplete 
gamma function, Kolmogorov differential equations, Kummer confluent hypergeometric function, skip-free Markov 
chain, species abundance and richness.} 

\tableofcontents

\section{Introduction}\label{sec-1}

This work is concerned with a variant of the \emph{birth and death chain} with constant intensities. 
Distinguished characteristic of the variant is that it allows negative jumps of any magnitude. It had already 
been analysed by Bruno de Finetti (see~\cite{deFinetti1957}) in a collectanea dedicated to the Italian 
mathematician Filippo Sibirani, that has had somewhat limited circulation and is now hard to find. This is 
one of the reasons why it might be worth revisiting that old work. Another reason, since de Finetti's paper 
is written in a rather informal style, is that it deserves a bringing into line with the modern theory of 
continuous-time Markov chains, that is de facto what the present one aims to do\footnote{In a correspondence 
(July 2025), Persi Diaconis kindly reminded me that the chain here dealt with represents a particular case of 
\emph{upward skip-free Markov chain}, a fact I omitted to mention out of ignorance. I have thus discovered a 
really interesting body of recent literature pertinent to the subject of the present paper. See, e.g., 
\cite{Choi2019}, \cite{Loeffen2025}. However, given its predetermined goal and the methodological approach 
adopted, I only have made minor changes to the first draft.}. But the most decisive factor in its 
revisitation is one sees a real possibility that it may be conducive to useful results concerning some 
aspects of the problem, arising in numerous scientific disciplines, of inferring the \emph{number} of 
undetected elements of some population on the basis of the detected elements of a proper part of it. From 
this point of view, in addition to the transparent, but somewhat limited example provided by de Finetti 
himself, one thinks of the topical issue of assessing the number of unseen specimens (abundance) of a given 
species (unseen species (richness), respectively) within a given area. Consider the case of a specific 
species of bird where detection occurs, e.g., by means of successive applications of the ``mark and 
recapture'' method. Then, the unknown number, say $X(t)$, of specimens which, beyond the ones marked during 
$[0,t]$, live in that area at time $t$, is an example of the number mentioned in the title. Its evaluation 
plays a crucial role in many ecological studies. Adaptation of the previous description to species richness 
is obvious. An example of a completely different nature is offered by the migrants who, during a certain 
interval time, say $[0,t]$, set sail from Africa headed for Italy along any Mediterranean route. Their number 
is undetectable, while there are data concerning arrivals, rescues at sea and shipwreck victims. Here, $X(t)$ 
is the number of migrants gone missing. In other words, it is the number of migrants beyond the ones who have 
been counted, during $[0,t]$, in one of the above three groups. Formalising the aforesaid phenomena according 
to de Finetti's model is tantamount to assuming that $X\coloneq \{X(t):t\ge 0\}$ forms a continuous-time 
Markov chain with state space $S=\Z^+\coloneq\{0,1,\ldots\}$, positive jumps of magnitude~$1$ -- 
corresponding to the entry of a new specimen, or species respectively, into the area, or a new embarkation -- 
and negative jumps of any magnitude -- corresponding to the number of specimens, or species, trapped and not 
yet marked, or migrants counted in any of the three groups as specified above. With respect to the extensive 
literature on the estimation of species abundance, or richness, the modeling described above, as well as 
possibile variants to correct its excessive simplicity, is decidedly unusual, at least to the present author 
knowledge. In fact, starting from the celebrated study by \cite{Fisher1943}, all approaches concentrate, even 
if from different viewpoints, on statistical modeling of procedures of observation, rather than the 
elementary factors in determining the number of undetected elements, and only at a later stage the 
statistical quantitative assessment of those factors. An analysis of the pros and cons is outside the scope 
of this paper. One limits oneself to giving a short list of recent papers which, besides specific original 
results, contain insightful accounts of the main existing different statistical approaches to estimating both 
richness and abundance of species: \cite{Orlitsky2016}; \cite{Baek2022}; \cite{Tekwa2023}, 
\cite{Mushagalusa2024}; \cite{Camerlenghi2024}. 

Returning to the present paper, one notes that the title is an almost literal translation of the Italian 
original, which is in agreement with the fact that it maintains the structure of the original. The main 
statements throughout the paper are labelled with capital letters in brackets, i.e., (\ref{th-A}), 
(\ref{th-B}),~\ldots\;. 

The rest of the paper is organised as follows. In Section~\ref{sec-2} one proves existence and uniqueness of 
the Markov chain corresponding to the infinitesimal description proposed by de Finetti. Section~\ref{sec-3} 
is devoted to the deduction of the corresponding transition function. Long-term properties (equilibrium, 
invariance, etc.)\ form the subject of Section~\ref{sec-4}. In particular, in Subsection~\ref{sub-4.1} one 
finds bounds on the error in approximating the equilibrium when that error is measured by means of both the 
Kolmogorov uniform distance and the Gini dissimilarity index. A bound on the absolute difference between 
moments of the same order is exhibited too. Section~\ref{sec-5} describes how to use the Markov chain under 
study in order to address some of the issues typically explored within the above-mentioned literature, 
including a statistical model which stems, in a consistent way, from the same chain. Section~\ref{sec-6} 
deals with the problem of estimating the unknown rates (parameters) of the infinitesimal description of the 
process, under the realistic assumption that no data are available except for jump times and magnitudes of 
negative jumps. The deduction of estimators from the above-mentioned statistical model, as it would be 
natural to do, requires an in-depth novel analysis which is outside the scope of the present work. Then, one 
confines oneself to considering the case in which jump times either remain unknown or are deliberately 
ignored, and, consequently, to adopting the conventional statistical model which sees the magnitudes as 
conditionally independent and identically distributed (i.i.d.)\ random numbers. Within such a model, one 
discusses definition and basic properties of a specific estimator obtained by the classical method of 
moments; see Subsections~\ref{sub-6.1} and \ref{sub-6.2}. The same estimator is reinterpreted, in 
Subsection~\ref{sub-6.4}, from a Bayesian standpoint. Section~\ref{sec-6} contains, here and there, remarks 
on both the numerical evaluation and approximation of the most relevant ``statistical'' functions which 
appear therein. The respective formal proofs are collected in a separate \nameref{sec-appendix}, split into a 
number of paragraphs referred to as \ref{app-A1}, \ref{app-A2}, \ldots\;. 

\section{Definition of the process}\label{sec-2}

As already mentioned, in \cite{deFinetti1957} it is assumed, without proof, that there exists a 
continuous-time Markov process with state-space $S=\Z^+$ and infinitesimal transition rates $q(i,j)$ given by 
\begin{equation}\label{eq-1}
  \begin{cases}
    q(i,i+1)=\lambda, q(i,i+k)=0 \quad \text{($i,k\in S$, $k\ge 2$)} \\
    q(i,i-k)=\mu \text{ or } 0 \quad\text{depending on $i\ge 1$ and $k=1,\ldots,i$, or $k>i\ge 0$} \\
    q(i,i)=-(\lambda+i\mu) \quad \text{($i\in S$)}
  \end{cases}
\end{equation} 
where $\lambda$ and $\mu$ are strictly positive numbers: $\lambda$, $\mu$ and $(\lambda+i\mu)$ represent the 
probability intensities to go from $i$ to $(i+1)$, from $i$ to $(i-k)$ and to leave $i$, respectively. For 
the sake of illustration, in the second example, $\lambda\,dt$ represents the first-order approximation to 
the probability that, during a length $dt$ time-interval, a new migrant boards, $\mu\,dt$ the same kind of 
approximation to the probability that any number of migrants be counted in one of the three classes 
(arrivals, rescue at sea, deaths at sea). 

The remaining part of the section focuses on the proof of the existence of the Markov process consistent with 
the infinitesimal description~\eqref{eq-1}. To this end, following the established literature on the subject 
(in particular, \cite{Norris1997} and \cite{Liggett2010}), one introduces

\begin{itemize}
  \item the set $\Omega$ of all right-continuous functions $\omega\colon \R^+\coloneq [0,+\infty)\to S$ 
      with \emph{finitely many jumps} in any \emph{bounded time-interval},
  \item the random function $X(t,\omega)\coloneq\omega(t)$ defined for every 
      $(t,\omega)\in\R^+\times\Omega$,
  \item its section $X(t)\coloneq X(t,\cdot\,)$ determined by $t$. 
\end{itemize} 

The set $S$ is endowed with the discrete topology and the corresponding $\sigma$-algebra, i.e.\ the power set 
$\mathcal{P}(S)$, while the $\sigma$-algebra $\mathcal{F}$ on $\Omega$ is the smallest such that the mapping 
$X(t)$ is measurable for every $t\in\R^+$. Moreover, for every $t\in\R^+$, $\mathcal{F}^{\,0}_t$ is defined 
to be the smallest $\sigma$-algebra w.r.t.\ which $X(s)$ is measurable for every $s\le t$, and then 
$\{\mathcal{F}_s \coloneq \bigcap_{t>s}\mathcal{F}^{\,0}_t:s\ge 0\}$ forms a right-continuous filtration and, 
of course, the process $\{X(t):t\ge 0\}$ turns out to be adapted to it. In this context, a \emph{Markov 
chain} on $S$ is meant as

\begin{enumerate}[leftmargin=*,label=(*)]
  \item \labeltext{*}{eq-*}a collection of probability measures $\{\mathbb{P}^x:x\in S\}$ on 
      $(\Omega,\mathcal{F})$ satisfying both  
      \[
        \mathbb{P}^x\{X(0)=x\}=1
      \]
      and the \emph{Markow property}
      \[
        \mathbb{P}^x(Y\circ X(t+s)\mid \mathcal{F}_s)=\mathbb{P}^{X(s)}(Y\circ X(t+s))\quad 
        \text{a.s.-$\mathbb{P}^x$}
      \]
      for all $x\in S$, all bounded measurable $Y\colon (\Omega,\mathcal{F})\to(\R,\mathcal{B}(\R))$, all 
      $s,t$ in $\R^+$.\footnote{Throughout the paper, one adopts the notation $\mu(f)\coloneq \int 
      f\,d\mu$.} 
\end{enumerate}

It is well-known that, as far as the construction of any Markov chain is concerned, the concept of transition 
function plays a key role. It is then worth recalling that any \emph{transition function} $p_t(x,y)$ is 
defined, for every $t\in\R^+$ and every $(x,y)\in S^2$, in such a way that
\begin{equation}\label{eq-**}
  p_t(x,y)\ge 0,\quad \sum_{y\in S}p_t(x,y)=1,\quad \lim_{t\downarrow 0^+}p_t(x,x)=p_0(x,x)=1 \tag{**}
\end{equation} 
hold true along with the Chapman-Kolmogorov equations
\begin{equation}\label{eq-***}
  p_{s+t}(x,y)=\sum_{z\in S} p_s(x,z)p_t(z,y). \tag{***}
\end{equation}

In point of fact, given the same measurable space $(\Omega,\mathcal{F})$ as above, along with a transition 
function $p_t(x,y)$, the products 
\[
  \prod_{j=0}^{n-1}p_{t_{j+1}-t_j}(x_j,x_{j+1})
\]
defined for every $0\eqcolon t_0<t_1<\ldots<t_n$, every $(x_1,\ldots,x_n)\in S^n$ and every 
$n\in\Z^+\setminus\{0\}$, form a system of consistent finite-dimensional distributions, whatever the initial 
state $x_0$ may be. A salient fact is that combination of the Kolmogorov extension theorem with 
right-continuity of the elements of $\Omega$ implies that there is one and only one probability measure 
$\mathbb{P}^x$ on $(\Omega,\mathcal{F})$ satisfying both the Markov property and 
\[
  \mathbb{P}^x\{X(t_1)=x_1,\ldots,X(t_n)=x_n\}=\prod_{j=0}^{n-1}p_{t_{j+1}-t_j}(x_j,x_{j+1})
\]
with $x_0=x$. This is nothing but a constructive definition of Markov chain that is rightly called 
\emph{transition probability definition}. Thus, returning to the original infinitesimal setting, it has to be 
checked whether the infinitesimal rates~\eqref{eq-1} determine a unique transition function or, to be more 
explicit, there is a unique transition function $p_t(x,y)$ that satisfies the \emph{Kolmogorov backwards 
equations}
\begin{equation}\label{eq-2}
  \begin{cases}
    \displaystyle{\frac{d}{dt}\, p_t(x,y)=\sum_{z\in S}q(x,z)p_t(z,y)\quad
    \text{($t>0$, $(x,y)\in S^2$)}} \\[1.2em]
    p_0(x,y)=\delta_{xy} \quad\text{($(x,y)\in S^2$)}
  \end{cases}
\end{equation} 
when the rates $q$ are the same as in~\eqref{eq-1}, and $\delta_{xy}$, as usual, stands for the Kronecker 
delta.

It is well-known that~\eqref{eq-2} has a minimal solution $p^*_t(x,y)$ which satisfies~\eqref{eq-***}, 
and~\eqref{eq-**} possibly except for stochasticity, i.e.\ the second relationship therein 
(\cite{Liggett2010}: Theorems~2.25, 2.26(a) and Proposition 2.21). In truth, these facts take place for any 
matrix $Q\coloneq\{ q(i,j):i,j\in S\}$ such that $q(i,j)\ge 0$ for $i\ne j$ and $\sum_j q(i,j)=0$. It is now 
interesting to wonder whether the minimal solution to~\eqref{eq-2}, when the elements of $Q$ are specified 
by~\eqref{eq-1}, is stochastic. Indeed, in such a case, $p^*_t(x,y)$ \emph{turns out to be the unique 
transition function satisfying}~\eqref{eq-2}, according to point~(b) in the aforementioned Theorem~2.26. To 
adequately discuss the issue, one introduces the \emph{jump matrix} $\Pi$ associated with the matrix $Q$ 
defined by~\eqref{eq-1}, that is 
\[
  \Pi\coloneq\biggl\{\pi(i,j)\coloneq\frac{q(i,j)}{c(i)}\text{ if $i\ne j$}, \pi(i,i)=0, c(i)\coloneq 
  -q(i,i)=\lambda+i\mu\biggr\}.
\]

In point of fact, the following statement is valid.

\medskip

\begin{regtheorem}\label{th-A}
Any discrete-time Markov chain on $S$, with initial distribution a strictly positive probability measure 
$\tau$ on $\mathcal{P}(S)$ and transition probabilities $\pi(i,j)$ given by the entries of $\Pi$, is 
irreducible and recurrent. Then, the minimal solution to~\eqref{eq-2}, $p^*_t(x,y)$, is stochastic. 
\end{regtheorem}

\begin{proof}[Proof of \textup{(\ref{th-A})}]
  For the validity of the latter statement see, e.g., \cite{Liggett2010}, Corollary~2.34(b). As to the former, 
  irreducibility is obvious while recurrence can be established by proving that the homogeneous system
  \begin{equation}\label{eq-3}
    \begin{cases}
      x_i=\sum_{j\ne 1}\pi(i,j)x_j & i\ne 1 \\[.5em]
      0\le x_i\le 1                & i\ne 1
    \end{cases}
  \end{equation}
  does not admit other solutions than the trivial one; see \cite{Billingsley1995}, Theorem~8.5. Now, it holds 
  that $x_0=0$ and, for any $(x_n)_{n\ge 2}$ satisfying the system, one obtains 
  $x_3-x_2=2x_2/\theta$, $x_{i+1}-x_i=(1+i/\theta)(x_i -x_{i-1})$ 
  for every $i\ge 3$, where $\theta=\lambda/\mu$. Thus, $x_n\uparrow +\infty$ if \mbox{$x_2>0$}, 
  which prevents it from being a solution to~\eqref{eq-3}, whilst $x_n=0$ for every $n=2,3,\ldots\,$, if 
  \mbox{$x_2=0$}.
\end{proof}

This paves the way for a direct application of the theorem of existence and uniqueness for Markov chains 
(see, e.g., \cite{Liggett2010}, Theorem~2.37). 

\medskip

\begin{regtheorem}\label{th-B}
There is a unique transition function $p^*_t(x,y)$ such that the elements $q$ given in~\eqref{eq-1} represent 
actual probability intensities, i.e.,
\[
  q(x,y)=\frac{d}{dt}\,p^*_t(x,y)\Big\lvert_{t=0} \qquad \textup{($x,y\in S^2$)}.
\]
Moreover, there is a unique Markov chain in the sense of \textup{(\ref{eq-*})} such that the measures 
$\mathbb{P}^x$ satisfy 
\[
  \mathbb{P}^x\{X(t)=y\}=p^*_t(x,y) \quad \textup{($x,y\in S^2$, $t\ge 0$)}.
\]
\end{regtheorem}

Proposition~(\ref{th-B}) confirms the substantial validity of the procedure followed in the original paper to 
define the process of interest. To complete the picture, here is a hint about the so-called \emph{jump 
chain}/\emph{holding time} description of the Markov chain under study. \emph{Jump times} are denoted by 
$J_0,J_1,\ldots$ and \emph{holding times} by $S_1,S_2,\ldots$\;: 
\[
  J_0\equiv 0, J_{n+1}\coloneq \inf\{t\ge J_n:X(t)\ne X(J_n)\} \qquad\text{($\inf\emptyset\coloneq +\infty$)}
\]
and
\[
  S_n\coloneq J_n-J_{n-1}
\]
for $n=1,2,\ldots$, with the proviso that $S_n\coloneq+\infty$ if $J_{n-1}=+\infty$, $X(+\infty)\coloneq 
X(J_n)$ if $J_{n+1}=+\infty$ for some $n$. It should be noted that, according to~(\ref{eq-*}), one has 
$\sup_{n}J_n=+\infty$, i.e, \emph{$X$ does not explode}. The discrete-time process $(Y_n)_{n\ge 0}$ defined 
by
\[
  Y_n\coloneq X(J_n)
\]
is called \emph{jump chain}. The relevant facts are:
\begin{itemize}
  \item $X(t)=Y_n$ for $J_n\le t < J_{n+1}$ ($n=0,1,\ldots)$.
  \item Conditional on $\{X(0)=x\}$, the jump chain $(Y_n)_{n\ge 0}$ is discrete-time Markov 
      $(\delta_x,\Pi)$, where $\delta_x$ stands for the point mass at $x$ and $\Pi$ is the same jump matrix 
      as in~(\ref{th-A}). 
  \item Conditional on $(Y_n)_{n\ge 0}$, the holding times $S_1, S_2,\ldots$ turn out to be independent 
      random numbers, and $S_n$ has exponential distribution of parameter $c(Y_{n-1})=\lambda+\mu Y_{n-1}$ 
      ($n=1,2,\ldots$). 
\end{itemize}

Conversely, given sequences $(Y_n)_{n\ge 0}$, $(S_n)_{n\ge 1}$ with $Y_0\equiv x$ and the same distributional 
properties as the above jumps and holding times, respectively, then the random function
\[
  t\mapsto Y_{\nu(t)} \quad\text{on $\{t:\nu(t)<+\infty\}$}
\]
where
\[
  \nu(t)\coloneq \begin{cases}
                   \min\{n\ge 0:S_1+\ldots+S_{n+1}>t\} & \text{if $\sum_{j\ge 1}S_j>t$} \\[.5em]
                   \infty & \text{otherwise}
                 \end{cases}
\]
forms a continuous-time Markov chain with the same transition function $p^*_t$ as in~(\ref{th-B}).

\section{Expression of the transition function}\label{sec-3}

Expressions for $p^*_t$ can be obtained by solving either \eqref{eq-2} or the \emph{Kolmogorov forward 
equations}
\begin{equation}\label{eq-4}
  \begin{cases}
    \dss{\frac{d}{dt}\,p_t(x,y)=\lambda p_t(x,y-1)-(\lambda+\mu y)p_t(x,y)+}\\[.5em] 
    \hspace{4cm}\dss{+\mu\sum_{\mathclap{z\ge y+1}}p_t(x,z)} \qquad 
    \text{($t>0$, $(x,y)\in S^2$, $p_t(x,-1)\equiv 0$)}  \\[2em]
    p_0(x,y)\coloneq \delta_{xy} \quad \text{($(x,y)\in S^2$)} 
  \end{cases}
\end{equation} 
since $p^*_t$ must be solution to \eqref{eq-4} in view of a well-known theorem (see, e.g. \cite{Norris1997}, 
Theorem~2.8.6). In fact, \cite{deFinetti1957} considers the equivalent form of~\eqref{eq-4} one obtains by 
summing the terms in the first line of \eqref{eq-4} w.r.t.\ $y$ over $\{0,1,\ldots,n-1\}$. In fact, putting 
\[
  R_t(\delta_x,n)\coloneq \sum_{y\ge n} p_t(x,y) \quad\text{if $t>0$}
\]
and
\[
  R_0(\delta_x,n)\coloneq \1\{n\le x\}
\]
for every $(x,n)$ in $S^2$, one obtains
\begin{equation}\label{eq-5}
  \begin{dcases}
    \frac{d}{dt}\,R_t(\delta_x,n)=\lambda R_t(\delta_x,n-1)-(\lambda+n\mu)R_t(\delta_x,n) \\[1em]
    n=1,2,\ldots
  \end{dcases}
\end{equation}
It is easy to see that the solution to \eqref{eq-5} given by formula \eqref{eq-3} in the original paper is 
misprinted: e.g., it violates the condition $R_t(\delta_x,n)\to R_0(\delta_x,n)$ as $t\to 0^+$. Then, new 
computations are made in the next subsection. 

\subsection{Equation of the generating function of the transition probability} 

The generating function 
\[
  g_t(x;z)\coloneq \sum_{n\ge 0} R_t(\delta_x,n)z^n \qquad\text{($t\ge 0$, $x\in S$, $-z_0\le z\le z_0$ for 
  some $z_0$ in $(0,1)$)}
\]
is used, via the theorem for the differentiation of power series, to transform \eqref{eq-5} into 
\[
  \sum_{n\ge 0}z^n\frac{d}{dt}\,R_t(\delta_x,n)=\lambda-\lambda(1-z)g_t(x;z)-\mu z\frac{\partial}{\partial z}\,
  g_t(x;z)
\]
which becomes
\begin{equation}\label{eq-6}
  \biggl(\frac{\partial}{\partial t}+\mu z\frac{\partial}{\partial z}\biggr)g_t(x;z)=\lambda-
  \lambda(1-z)g_t(x;z)
\end{equation}
by force of the following argument. Letting $f_t(n;z)\coloneq z^n R_t(\delta_x,n)$, one can write
\[
  g_t(x;z)\coloneq \sum_{n\ge 0} R_t(\delta_x,n)z^n=\int_{S} f_t(n;z)\,\nu(dn)
\]
where $\nu$ stands for the counting measure on $S$, and
\[
  \abs[\bigg]{\frac{\partial}{\partial t}f_t(n;z)}=
  \abs{z}^n\abs[\bigg]{\frac{\partial}{\partial t}R_t(\delta_x,n)}\le
  (\lambda+n\mu)\cdot\abs{z_0}^n.
\]
Since the majorant function is independent of $t$, and
\[
  \int_S (\lambda+n\mu)\abs{z_0}^n\,\nu(dn)<+\infty
\]
then differentiation under integral sign
\[
  \frac{\partial}{\partial t}\int_{S}f_t(n;z)\,\nu(dn)=\int_S\frac{\partial}{\partial t}f_t(n;z)\,\nu(dn)
\]
is valid for $t>0$ and $-z_0\le z\le z_0$; see \cite{Billingsley1995}, Theorem~16.8~(ii). Now, one has to 
find solutions to~\eqref{eq-6} which satisfy the initial condition 
\begin{equation}\label{eq-7}
  \dss{g_0(x;z) = \sum_{n\ge 0} \1\{n\le x\}z^n = \sum_{n=0}^{x}z^n} \qquad\text{($\abs{z}<z_0$)}.
\end{equation} 
Equation \eqref{eq-6} is equivalent to the system of ordinary differential equations
\[
  dt=\frac{dz}{\mu z}=\frac{dg_t}{\lambda(1-(1-z)g_t)}.
\]
Integration of
\[
  dt=\frac{dz}{\mu z}
\]
gives
\begin{equation}\label{eq-8}
  c_1=ze^{-\mu t}.
\end{equation}
Taking
\[
  \frac{dz}{\mu z}=\frac{dg_t}{\lambda(1-(1-z)g_t)}
\]
as second equation of the system, i.e.
\[
  \frac{dg_t}{dz}=-\theta\,\frac{1-z}{z}g_t+\frac{\theta}{z}\qquad \text{($\theta=\lambda/\mu$)}
\]
one has
\begin{align*}
  g_t(x;z) & = \abs{z}^{-\theta} e^{\theta z}\biggl[c_2+\theta\int_{0}^{z}\frac{1}{x}\abs{x}^\theta e^{-\theta x}\,dx\biggr] \\[.3em]
           & = \abs{z}^{-\theta} e^{\theta z}\biggl[c_2+\theta\abs{z}^\theta\int_{0}^{1}u^{\theta-1} e^{-\theta zu}\,du\biggr] \\[.3em]
           & = \abs{z}^{-\theta} e^{\theta z}\biggl[c_2+\abs{z}^\theta e^{-\theta z}\Phi(1,\theta+1,\theta z)\biggr]
\end{align*}
where $\Phi$ denotes the \emph{Kummer confluent hypergeometric function}. Hence, letting $t=0$, which entails 
$z=c_1$ in view of \eqref{eq-8}, and recalling the initial condition \eqref{eq-7}, one obtains
\[
  c_2=\abs{c_1}^\theta e^{-\theta c_1}g_0(x;c_1)-\abs{c_1}^\theta e^{-\theta c_1}\Phi(1,\theta+1,\theta c_1)
\]
where
\[
  g_0(x;c_1)=\sum_{j=0}^{x} c_1^j.
\]
Then, since $c_1=ze^{-\mu t}$ (from~\eqref{eq-8}), it turns out that
\[
  c_2=\abs{z}^\theta e^{-\theta\mu t-\theta ze^{-\mu t}}\biggl[\,\sum_{j=0}^{x}z^je^{-\mu tj}-\Phi(1,\theta+1,\theta ze^{-\mu t})\biggr]
\]
and, consequently,
\[
  g_t(x;z)=\Phi(1,\theta+1,\theta z)+e^{-\theta z(e^{-\mu t}-1)-\theta\mu t}\biggl[\,\sum_{j=0}^{x}z^je^{-\mu tj}-
  \Phi(1,\theta+1,\theta ze^{-\mu t})\biggr]
\]
which, in view of the power series definition of the Kummer function $\Phi$, becomes
\[
  g_t(x;z)=\sum_{n\ge 0} z^n R^*_t(\delta_x,n)
\]
with
\begin{equation}\label{eq-9}\begin{split}
  \hspace{3em}R^*_t(\delta_x,n) & =\frac{\theta^n}{(\theta+1)_n}+ \\[.5em]
   & +e^{-(\theta+n)\mu t}\sum_{\rho=0}^{n}\frac{(\theta(e^{\mu t}-1))^\rho}{\rho!} 
   \biggl(\1_{[n-\rho,+\infty)}(x)-\frac{\theta^{n-\rho}}{(\theta+1)_{n-\rho}}\biggr)
\end{split}
\end{equation}
$(\theta)_n$ being the Pochhammer symbol for the \emph{rising factorial}, i.e., $(\theta)_n\coloneq 
\theta(\theta+1)\ldots (\theta+n-1)$ if $n\ge 1$, and $(\theta)_0\coloneq 1$. Then,

\begin{regtheorem}\label{th-C}
  The transition function of the Markov chain, whose existence has been
   established in~\textup{(\ref{th-B})}, is given by
  \[
    p^*_t(x,y)\equiv R^*_t(\delta_x,y)-R^*_t(\delta_x,y+1)
  \]
  for every $x,y$ in $S$ and $t\ge 0$, where $R^*_t$ is defined by~\eqref{eq-9}.
\end{regtheorem}

In the rest of the paper, given any probability $\tau$ on $S$, the symbols $(\tau,\{\bfP^x:x\in S\})$ and 
$\bfP^{(\tau)}$ will designate the Markov chain $X$ with initial distribution $\tau$ and transition function 
$p^*_t(x,y)$, and the corresponding probability measure on $(\Omega,\mathcal{F})$, respectively. Then,
\begin{align*}
  R^*_t(\tau,n) & \coloneq \bfP^{(\tau)}\{X(t)\ge n\}=\sum_{x\ge 0}\tau(x)R^*_t(\delta_x,n) 
  \qquad\text{($n\in S$)}  \\[.5em]
       & = I(n)+e^{-(\theta+n)\mu t}\sum_{\rho=0}^{n}\frac{(\theta(e^{\mu t}-1))^\rho}{\rho!}\Delta(n-\rho)
\end{align*}
where
\[
  I(n)\coloneq \frac{\theta^n}{(\theta+1)_n} \qquad \text{($n\in S$)}
\]
and
\[
  \Delta(n)\coloneq\tau([n,+\infty))-I(n) \qquad \text{($n\in S$)}.
\]

An immediate consequence of these definitions is that
\begin{equation}\label{eq-10}
  \pi^*(n)\coloneq I(n)-I(n+1)=\frac{\theta^n}{(\theta+1)_{n+1}}(n+1) \qquad \text{($n\in S$)}
\end{equation}
is a probability density function w.r.t.\ the counting measure on $(S,\mathcal{P}(S))$. The corresponding 
probability measure, denoted by the same symbol, has \emph{finite moments of any order}. Indeed, for every 
$\rho>0$, one has 
\begin{equation}\label{eq-11}
  \sum_{n\ge 1} n^\rho\pi^*(n)=\rho\sum_{n\ge 1}\int_{n-1}^{n}x^{\rho-1}I(n)\,dx =\sum_{n\ge 1}I(n)\{n^\rho -
  (n-1)^\rho\}.
\end{equation}
Since
\[
  0<I(n)\{n^\rho-(n-1)^\rho\}\le n^\rho\,\frac{\theta^n}{n!} \qquad \text{($n\ge 1$)}
\]
then the series has a finite sum. 

Moreover, setting $z\coloneq \theta(e^{\mu t}-1)$, one gets 
\begin{multline*}
  \sum_{n\ge 1}\{n^\rho-(n-1)^\rho\}
  \abs[\bigg]{e^{-(\theta+n)\mu t}\sum_{\rho=0}^{n}\frac{z^\rho}{\rho!}\Delta(n-\rho)} \le \\[.5em]
  \le \sum_{n\ge 1}\{n^\rho-(n-1)^\rho\}e^{-(\theta+n)\mu t+z}<+\infty 
\end{multline*}
which, combined with the previous result, entails
\begin{equation}\label{eq-12}
  \bfP^{(\tau)} (\abs{X(t)}^\rho)<+\infty
\end{equation}
\emph{for every starting distribution $\tau$, every $\rho>0$ and $t>0$}. Of course, \eqref{eq-12} extends to 
$t=0$ if $\tau$ has finite $\rho$-th moment.

Other properties of $(\tau,\{\bfP^x:x\in S\})$, in particular long-time properties, are described in the next 
section. 

\section{Long-time properties of the process}\label{sec-4}

It is plain to see that $\pi^*$ is the equilibrium distribution, i.e., the limiting probability distribution 
of $p^*_t(x,\cdot\,)$, as $t\to+\infty$, for all $x$. Some of its more interesting properties are described 
in the following two propositions. The former concerns stationarity. 

\begin{regtheorem}\label{th-D}
  The limiting probability distribution $\pi^*$ is invariant for the Markov chain with transition 
  function $p^*_t$, that is
  \[
    \pi^*(y)=\sum_{x\ge 0}\pi^*(x)p^*_t(x,y) \qquad \textup{($y\in S$, $t\ge 0$)}.
  \]
  Moreover, $(\pi^*,\{\bfP^x:x\in S\})$ is a stationary Markov chain, i.e., the joint distribution
  \[
    \sum_{x\in S}\pi^*(x)\,\bfP^x\{X(t_1+s)\in A_1,\ldots,X(t_n+s)\in A_n\} \qquad 
    \textup{($A_j\subset S$, $j=1,\ldots,n$)}
  \]
  is independent of $s$ whenever $0\le t_1+s<t_2+s<\ldots<t_n+s$.
\end{regtheorem} 

The latter proposition pertains to the class structure of the chain and ergodicity.

\begin{regtheorem}\label{th-E}
  Every state $x\in S$ is positive recurrent, i.e.,
  \[
    \bfP^x\bigl(\textup{$\{t\ge 0:X(t)=x\}$ is unbounded}\bigr)=1
  \]
  and the expectation $\bfP^x$ of the return time
  \[
    T_x\coloneq \inf\{t>J_1:X(t)=x\}
  \]
  obeys
  \[
    \bfP^x(T_x)=\frac{(\theta+1)_x}{\mu(\theta+x)(1+x)\theta^x}.
  \]
  Moreover, the ergodic condition holds true, i.e.,
  \[
    \bfP^{(\tau)}\biggl\{\frac{1}{t}\int_{0}^{t}f(X(s))\,ds\to \sum_{x\ge 0}\pi^*(x)f(x) 
    \textup{ as $t\to+\infty$}\biggr\}=1
  \]
  for any $\tau$ and every bounded function $f\colon S\to\R$.
\end{regtheorem}

Propositions (\ref{th-D}) and (\ref{th-E}) follow from plain combinations of (\ref{th-A}) with some basic 
results, for which the reader is referred to \cite{Norris1997}: Theorems 3.5.1, 3.5.2, 3.5.3, 3.8.1, 
\cite{Liggett2010}, Subsection~2.6.1. 

\subsection{Bounds on the error in approximating the equilibrium}\label{sub-4.1}

In this subsection the speed of approach to equilibrium is considered from two different points of view: 
convergence of the cumulative distribution of $X(t)$ and convergence of its moments. The corresponding 
results also provide expressive bounds on the error in approximating for fixed $t$.

\begin{sloppypar}
The first statement is about the Kolmogorov \emph{uniform distance} between $\Pi^*(x) \coloneq \sum_{n\le 
x}\pi^*(n)$, $x\in\R$, and the cumulative distribution function $F^{(\tau)}_{X(t)}$ of $X(t)$, i.e., 
$F^{(\tau)}_{X(t)}(x)\coloneq \bfP^{(\tau)}\{X(t)\le x\}$. 
\end{sloppypar}

\medskip

\begin{regtheorem}\label{th-F}
  If $T(x)\coloneq\sum_{n\le x}\tau(n)$, $x\in\R$, then
  \begin{multline*}
    \sup\{\abs{\Pi^*(x)-F^{(\tau)}_{X(t)}(x)}:x\in\R\}\\[.5em]
     \le\sup\{\abs{\Pi^*(x)-T(x)}:x\in\R\} \exp\{-\mu t-\theta(e^{\mu t}+\mu t-1)\}.
  \end{multline*}
\end{regtheorem}

\begin{proof}[Proof of \textup{(\ref{th-F})}]
  The uniform distance between $\Pi^*$ and $F^{(\tau)}_{X(t)}$ is equal to
  \begin{align*}
    \sup_{n\ge 1}\abs{R^*_t(\tau,n)-I(n)} & \le\sup_n 
    \Delta(n)\cdot\sup_n e^{-(\theta+n)\mu t}e_{n-1}(\theta(e^{\mu t}-1)) 
    \shortintertext{(from the end of Section~\ref{sec-3}, with the proviso that $e_{\nu}(x)\coloneq 
    \sum_{k=0}^{\nu}x^k/k!$)} \\
     & =\sup_n\Delta(n) \cdot
     \sup_n e^{-\theta\mu t}\sum_{\rho=0}^{n-1}e^{-\mu t(n-\rho)}\frac{(\theta(1-e^{-\mu t}))^\rho}{\rho!} \\[1em]
     & \le \sup_n \Delta(n)e^{-\theta\mu t-\mu t+\theta(1-e^{-\mu t})}. \qedhere
  \end{align*}
\end{proof}

The next result concerns the error of the $m$-th moment of $\Pi^*$ in approximating the homologous moment of 
$F^{(\tau)}_{X(t)}$.

\begin{regtheorem}\label{th-G}
  If the starting distribution $\tau$ has finite $m$-th moment, for some integer $m\ge 1$, then
  \[
    \abs[\bigg]{\int_\R x^m\,d\bigl(F^{(\tau)}_{X(t)}(x)-\Pi^*(x)\bigr)}\le m\,\mathcal{K}_{h_m}(\tau,\pi^*)\cdot 
    \exp\{-\mu t-\theta(e^{-\mu t}+\mu t-1)\}\qquad \textup{($t\ge 0$)}
  \]
  where $h_m$ is a real-valued function and $\mathcal{K}_{h_m}(\tau,\pi^*)$ a real number defined by
  \[
    h_m(x)\coloneq\sum_{\rho\ge 0}\frac{e^{-\theta}}{\rho!}\theta^\rho(x+\rho)^{m-1} \qquad \textup{($x\ge 0$)}
  \]
  and
  \[
    \mathcal{K}_{h_m}(\tau,\pi^*)\coloneq \int_{0}^{+\infty}h_m(x)\abs{T(x)-\Pi^*(x)}\,dx
  \]
  respectively.
\end{regtheorem} 

\begin{remark}
$\mathcal{K}_{h_m}$ is an example of \emph{Kantorovich-Rubinstein functional}. It represents the value of the 
minimal (total) translocation cost  in the case of transits permitted (masses~$q$ on $S^2$ such that 
$q(\{n\}\times S)-q(S\times\{n\})=\tau(n)-\pi^*(n)$, for every $n\in S$), when the unit cost $c$ of 
transportation from $x$ to $y$ is given by 
\[
  c(x,y)=\abs{x-y}\max(h_m(x),h_m(y)).
\]

For a comprehensive account of these concepts, the reader is referred to \cite{Rachev2013}, Chapter~5. 
\end{remark}

\begin{proof}[Proof of \textup{(\ref{th-G})}]
  Recalling that the $m$-th moment of a probability distribution $p$ \emph{supported by some subset of $\R^+$}, 
  with cumulative distribution function $F_p$, is presentable as
  \[
    \int_{\R} x^m\, p(dx)=m\int_{0}^{+\infty}x^{m-1}[1-F_p(x)]\, dx
  \]
  then, by force of \eqref{eq-11}-\eqref{eq-12}, the absolute difference between moments to be investigated 
  can be written and majorised as follows
  \begin{align*}
    \abs[\bigg]{\int_\R x^m\,d(F^{(\tau)}_{X(t)}-\Pi^*)}=
    \abs[\bigg]{\sum_{k\ge 0}\bigl((k+1)^m-k^m\bigr)\big(R^*_t(\tau,k+1)-I(k+1)\big)} \\[.5em]
    \le f(t;\theta,\mu)\sum_{\rho\ge 0}\mathcal{P}_{g(t)}(\rho)\sum_{\nu\ge 1}[(\rho+\nu)^m-(\rho+\nu-1)^m]
    \cdot \abs{T(\nu-1)-\Pi^*(\nu-1)}
  \end{align*}
\begin{sloppypar}
\noindent (where $\mathcal{P}_{g(t)}$ denotes the Poisson distribution with parameter $g(t)\coloneq 
\theta(1-e^{-\mu t})$, $f(t;\theta,\mu)\coloneq\exp\{g(t)-t(\lambda+\mu)\}$) 
\end{sloppypar}
\begin{align*}
    & = f(t;\theta,\mu)\sum_{\rho\ge 0}\mathcal{P}_{g(t)}(\rho)\sum_{\mu\ge 1}m\int_{\nu-1}^{\nu}(\rho+x)^{m-1}\,dx 
     \cdot \abs{T(\nu-1)-\Pi^*(\nu-1)} \\[.5em]
    &  =f(t;\theta,\mu)m D_m(t)
\end{align*}
where $D_m(t)\coloneq\int_{0}^{+\infty}\sum_{\rho\ge 0}\mathcal{P}_{g(t)}(\rho) 
(\rho+x)^{m-1}\cdot\abs{T(x)-\Pi^*(x)}\,dx$. The last equality is a consequence of the Fubini theorem applied 
to the integral of the positive function $(\rho,x)\mapsto (\rho+x)^{m-1}\abs{T(x)-\Pi^*(x)}$, $(\rho,x)\in 
S\times\R^+$, w.r.t.\ the measure $\mathcal{P}_{g(t)}\times\mathrm{Leb}$. To complete the proof it suffices 
to check the validity of the following relations, in which $S(\cdot\,,\cdot\,)$ denotes the generic Stirling 
number of the second kind. 
\begin{align*}
  G(t,x) & \coloneq \sum_{\rho\ge 0}\mathcal{P}_{g(t)}(\rho)(\rho+x)^{m-1} 
  \qquad \text{(finite for every $t\ge 0$, $x\ge 0$)} \\[.5em]
         & = \sum_{j=0}^{m-1} \binom{m-1}{j} x^{m-1-j}\sum_{k=0}^{j}S(j,k)g(t)^k \\[.5em]
         & \le \sum_{j=0}^{m-1} \binom{m-1}{j} x^{m-1-j}\sum_{k=0}^{j}S(j,k)\theta^k \\[.5em]
         & = \sum_{\rho\ge 0}\mathcal{P}_\theta(\rho)(\rho+x)^{m-1}\eqcolon h_m(x).
\end{align*} 
Then,
\begin{align*}
  0\le D_m(t) & =\int_{0}^{+\infty}G(t,x)\abs{T(x)-\Pi^*(x)}\,dx \\[.5em]
              & \le \int_{0}^{+\infty} h_m(x)\abs{T(x)-\Pi^*(x)}\, dx = \mathcal{K}_{h_m}(\tau,\pi^*) \\[.5em]
    & = \sum_{j=0}^{m-1}\sum_{k=0}^{j} \binom{m-1}{j} S(j,k)\,\theta^k\int_{0}^{+\infty}x^{n-1-j} 
    \abs{T(x)-\Pi^*(x)}\,dx \notag \\[.5em]
    & < +\infty
\end{align*}
since, by hypothesis,
\[
  \int_{0}^{+\infty} x^{m-1-j} \abs{T(x)-\Pi^*(x)}\,dx <+\infty
\]
for $j=0,\ldots,m-1$.
\end{proof}

\begin{remark}
If $m=1$ and $\bfP^{(\tau)}_{X(t)}$ denotes the probability distribution (under $\bfP^{(\tau)}$) of $X(t)$, 
then $\mathcal{K}_{h_m}(\bfP^{(\tau)}_{X(t)},\pi^*)$ reduces to the \emph{Gini dissimilarity index} 
\[
  \mathcal{G}(\bfP^{(\tau)}_{X(t)},\pi^*)=\int_\R\abs{F^{(\tau)}_{X(t)}(x)-\Pi^*(x)}\, dx
\]
that, historically, can be viewed as a prototype of Kantorovich functional. See, once again, 
\cite{Rachev2013}. Then, the same argument as in the previous proof can be used to validate the inequality 
\[
  \mathcal{G}(\bfP^{(\tau)}_{X(t)},\pi^*)\le\exp\{-\mu t-\theta(e^{-\mu t}+\mu t-1)\}\mathcal{G}(\tau,\pi^*)
\]
\emph{provided that $\tau$ has finite expectation}.
\end{remark}

\section{Outline of potential practical uses of the model}\label{sec-5}

Here is a description of how the Markov chain analysed in the last three sections might be used to tackle 
some of the problems outlined in the Introduction. One notes, preliminarily, that the following treatment 
rests on the realistic assumption that \emph{times and magnitudes of negative jumps are the sole facts one 
can actually observe}. Then -- denoted by $\bft^{(n)}$ and $\bfd^{(n)}$ the vectors $(t_0,\ldots,t_n)$ and 
$(d_1,\ldots,d_n)$, respectively, where $0=t_0<\ldots<t_n$, $d_k\in S\setminus\{0\}$ for $k=1,\ldots,n$, and 
$n\ge 1$ -- one considers the events 
\begin{align*}
  H_{\bft^{(n)}} & \coloneq \bigcap_{k=1}^n\{\text{a \emph{negative} jump occurs at $t_k$}\} \\[.5em]
  H_{\bft^{(n)},\bfd^{(n)}} & \coloneq \bigcap_{k=1}^n\{\text{a \emph{negative} jump of magnitude 
                                $d_k$ occurs at $t_k$}\}
\end{align*}
It is clear that the accomplishment of the aim set at the beginning of this section largely depends on the 
analysis of the (conditional) probabilities $\bfP^{(\tau)}(\bigcap_{j=1}^n\{X(t_j)=x_j\mid 
H_{\bft^{(n)},\bfd^{(n)}})$ and $\bfP^{(\tau)}(\bigcap_{j=1}^n\{\abs{\Delta X(t_j)}=d_j,X(t_j)=x_j\}\mid 
H_{\bft^{(n)}})$, where, assuming that $X(t-0)\ne X(t)$, $\Delta X(t)\coloneq X(t-0)-X(t)$ and $\abs{\Delta 
X(t)}$ stand for size and magnitude, respectively, of the jump at $t>0$. Now, by simple reasoning, one can 
prove

\begin{regtheorem}\label{th-H}
Let $x_0\coloneq x$ and, $\forall n\in\Z^+\setminus\{0\}$, let $(x_1,\ldots,x_n)$ belong to $S^n$, 
$\bfd^{(n)}\coloneq (d_1,\ldots,d_n)$ to $(S\setminus\{0\})^n$, $\bft^{(n)}\coloneq (t_0,\ldots,t_n)$ with 
$0=t_0<\ldots<t_n$, and $D^{(\tau)}_n(\bft^{(n)},\bfd^{(n)})\coloneq \sum_{x\in S}\tau(x)\sum_{x_1\ge 0} 
\ldots\sum_{x_n\ge 0}p^*_{t_1}(x,x_1+d_1)\ldots p^*_{t_n-t_{n-1}}(x_{n-1},x_n+d_n)$. Then, 
\begin{multline}\label{eq-13}
\bfP^{(\tau)}(\abs{\Delta X(t_1)}=d_1,X(t_1)=x_1,\ldots,\abs{\Delta X(t_n)}=d_n,X(t_n)=x_n\mid H_{\bft^{(n)}})= 
\\[.5em]
=\frac{1}{\sum_{d_1\ge 1}\ldots\sum_{d_n\ge 1}D^{(\tau)}_n(\bft^{(n)},\bfd^{(n)})}\sum_{x\ge 0}\tau(x) 
p^*_{t_1}(x,x_1+d_1)\cdot\ldots \\[.5em]
\ldots\cdot p^*_{t_n-t_{n-1}}(x_{n-1},x_n+d_n);
\end{multline}
\begin{multline}\label{eq-14}
\bfP^{(\tau)}(\abs{\Delta X(t_1)}=d_1,\ldots,\abs{\Delta X(t_n)}=d_n\mid H_{\bft^{(n)}})= \\[.5em]
=\frac{1}{\sum_{d_1\ge 1}\ldots\sum_{d_n\ge 1}D^{(\tau)}_n(\bft^{(n)},\bfd^{(n)})}
D^{(\tau)}(\bft^{(n)},\bfd^{(n)});
\end{multline}
\begin{multline}\label{eq-15}
\bfP^{(\tau)}(X(t_1)=x_1,\ldots,X(t_n)=x_n\mid H_{\bft^{(n)},\bfd^{(n)}})= \\[.5em]
=\frac{1}{D^{(\tau)}_n(\bft^{(n)},\bfd^{(n)})}\sum_{x\ge 0}\tau(x) 
p^*_{t_1}(x,x_1+d_1)\cdot\ldots\cdot p^*_{t_n-t_{n-1}}(x_{n-1},x_n+d_n).
\end{multline}
\end{regtheorem} 

Distribution \eqref{eq-14} may play a role as a basis for statistical inference, e.g., estimation of the 
generally unknown parameters $\theta$ and $\mu$ present in the expression of the transition function $p^*_t$. 
Indeed, consistently with the basic initial assumption about what is actually observable, for each 
$\bfd^{(n)}$ and $\bft^{(n)}$, \eqref{eq-14}, considered as a function of $(\theta,\mu)$, represents the 
\emph{likelihood function}, conditionally on the hypothesis that observations (negative jumps) occur at times 
$t_1,\ldots,t_n$. However, the rather laborious development of such a statistical idea goes beyond the aim of 
the present paper. A decidely more conventional likelihood function will be considered, along with some of 
its consequences, in Section~\ref{sec-6}, where one imagines situations in which jumping times are 
unavailable data.

Distribution \eqref{eq-15} can be used, e.g., to evaluate the probability that, immediately after the last 
observation of new specimens of a given species, the number of unseen elements of that very same species is 
not less than an arbitrarily fixed integer $\xi$, conditionally on the fact that those observations occurred 
at times $t_1,\ldots,t_n$ and consisted in $d_1,\ldots,d_n$ new specimens, respectively: 
\begin{multline}\label{eq-16}
  \bfP^{(\tau)}(X(t_n)\ge\xi\mid H_{\bft^{(n)},\bfd^{(n)}})=\frac{1}{D^{(\tau)}(\bft^{(n)},\bfd^{(n)})}
  \sum_{x\ge 0}\ldots\sum_{x_{n-1}\ge 0}\tau(x)p^*_{\Delta t_1}(x,x_1+d_1)\ldots \\[.5em]
  \ldots p^*_{\Delta t_{n-1}}(x_{n-2},x_{n-1}+d_{n-1})R^*_{\Delta t_n}(\delta_{x_{n-1}},d_n+\xi) \\[.5em]
  =\sum_{x_{n-1}\ge 0}m_{n-1}(x^{(\tau)}_{n-1};\bfd^{(n)},\Delta\bft^{(n)})\bfP^{x_{n-1}}(X(\Delta t_n)\ge 
  d_n+\xi\mid X(\Delta t_n)\ge d_n)
\end{multline} 
where $\Delta t_j\coloneq t_j-t_{j-1}$ ($j=1,\ldots,n$), $\Delta\bft^{(n)}=(\Delta t_1,\ldots,\Delta t_n)$,
\begin{multline*}
  m^{(\tau)}_{n-1}(x_{n-1};\bfd^{(n)},\Delta\bft^{(n)})\coloneq\frac{1}{D^{(\tau)}(\bft^{(n)},\bfd^{(n)})}
  \sum_{x\ge 0}\ldots\sum_{x_{n-2}\ge 0}\tau(x) \\[.5em]
  \ldots p^*_{\Delta t_{n-1}}(x_{n-2},x_{n-1}+d_{n-1})R^*_{\Delta t_n}(\delta_{x_{n-1}},d_n) \qquad
  \text{($n=2,1,\ldots$)}
\end{multline*}
\[
  m^{(\tau)}_0(x;d_1,t_1)\coloneq \frac{\tau(x)R^*_{t_1}(\delta_x,d_1)}{D^{(\tau)}(t_1,d_1)}
\]
and, of course,
\[
  \bfP^{x_{n-1}}(X(\Delta t_n)\ge d_n+\xi\mid X(\Delta t_n)\ge d_n)=
  \frac{R^*_{\Delta t_n}(\delta_{x_{n-1}}, d_n+\xi)}{R^*_{\Delta t_n}(\delta_{x_{n-1}},d_n)} \qquad 
  \text{($n=1,2,\ldots$)}.
\]
It is worth noticing that equality
\[
  \sum_{x\ge 0}m^{(\tau)}_{n-1}(x;\bfd^{(n)},\bft^{(n)})=1
\]
holds true for every $n=1,2,\ldots$\;, and $m^{(\tau)}_{n-1}(x;\bfd^{(n)},\bft^{(n)})>0$ for every $x,n$. 
Further useful information is also provided by the \emph{expectation} 
\begin{equation}\label{eq-17}
  \bfP^{(\tau)}(X(t_n)\mid H_{\bft^{(n)},\bfd^{(n)}})=\sum_{s\ge 0}m^{(\tau)}_{n-1}
  (s;\bfd^{(n)},\Delta\bft^{(n)})\sum_{k\ge 1}
  \frac{R^*_{\Delta t_n}(\delta_s, k\vee d_n)}{R^*_{\Delta t_n}(\delta_s,d_n)}-d_n.
\end{equation}

Reinterpretation of the previous probability evaluations would be perfectly justifiable from perspectives 
quite different than the species abundance, like, for example, the number of unseen migrants along the 
Mediterranean routes during a certain period. In any case, they represent probabilities of events 
\emph{conditional} on either the event $H_{\bft^{(n)}}$ or the event $H_{\bft^{(n)},\bfd^{(n)}}$, evaluated 
for every $\bft^{(n)}$ ($(\bft^{(n)},\bfd^{(n)})$, respectively) and for every $n$. One will alternately 
confine oneself to considering that distribution which corresponds to the specific $\bft^{(n)}$ 
($(\bft^{(n)},\bfd^{(n)})$, respectively) which turns out to be true. But it is quite natural to think of 
questions which require the evaluation of ``absolute'' probabilities like when, e.g., one is wondering about 
the probability distribution of the first jump time, jointly with its magnitude and the number of unseen 
elements immediately after its occurrence. As far as the calculus of these probabilities is concerned, the 
jump chain/holding times description, recalled at the end of Section~\ref{sec-2}, can provide proper tools. 
In the same notation introduced therein, completed by $\bar{J}_k$ to denote the $k$-th \emph{negative} jump 
time ($k=1,2,\ldots$), fix $\nu$ in $\Z^+$ and assume that $\bar{J}_1$ corresponds to the $(\nu+1)$-th jump 
of the chain. The $\bfP^x$-probability that such an event occurs jointly with $\{\abs{\Delta 
X(\bar{J}_1)}=d\}$ is given by 
\begin{equation}\label{eq-18}
  \frac{\lambda^\nu\mu}{(\lambda+\mu x)\ldots(\lambda+\mu(x+\nu))}\1\{1\le d\le x+\nu\}
  =\frac{\theta^\nu}{(\theta+x)_{\nu+1}}\1\{1\le d\le x+\nu\}
\end{equation} 
which can also be reinterpreted as $\bfP^x\{X(\bar{J}_1)=x+\nu-d,\abs{\Delta X(\bar{J}_1)}=d\}$. If 
$\bar{J}_1$ corresponds to the $(\nu+1)$-th jump of the chain, then its probability distribution is the same 
as that of the sum of $(\nu+1)$-independent exponentially distributed random numbers of parameter 
$\lambda+\mu(x+j)$, $j=0,\ldots,\nu$, respectively. Then, by the convolution theorem, the \emph{Laplace 
transform} of the density of that sum is
\begin{align*}
  p\mapsto & \prod_{j=0}^{\nu}\frac{\lambda+\mu(x+j)}{p+\lambda+\mu(x+j)}\qquad\text{($p>-\lambda$)} \\[.5em]
    & =\Gamma\begin{bmatrix}
            \theta+x+\nu+1 \\
            \theta+x
             \end{bmatrix}\cdot
       \Gamma\begin{bmatrix}
               \dfrac{p}{\mu}+\theta+x \\[1em]
               \dfrac{p}{\mu}+\theta+x+\nu+1 
             \end{bmatrix}
\end{align*} 
whose determining function is given by
\begin{equation}\label{eq-19}
  \R\ni t\mapsto\frac{\mu}{B(\nu+1,\theta+x)}e^{-(\theta+x)\mu t}(1-e^{-\mu t})^\nu\1\{t>0\}
\end{equation}
where $B$ stands for the \emph{Euler Beta function}. Then, combination of \eqref{eq-19} with \eqref{eq-18} 
gives
\begin{multline*}
  \bfP^x\{\bar{J}_1>t_1,X(\bar{J}_1)=x+\nu-d,\abs{\Delta X(\bar{J}_1)}=d\}= \\[.5em]
  =\frac{\theta^\nu}{(\theta+x)_{\nu+1}}\int_{t_1\vee 0}^{+\infty}\frac{\mu}{B(\nu+1,\theta+x)}
  e^{-(\theta+x)\mu t}(1-e^{-\mu t})^\nu\, dt\cdot\1\{1\le d\le x+\nu\}
\end{multline*} 
for every $x\in S$, $\nu\in\Z^+$, $d\in\Z^+\setminus\{0\}$ and $t_1\in\R$, i.e,
\begin{regtheorem}\label{th-I}
Given any $x\in S$, the function $f$ defined on $\R^+\times S\times (S\setminus\{0\})$ by
\begin{equation}\label{eq-20}
  f(x;t,x_1,d)\coloneq\frac{\mu\theta^{x_1+d-x}}{(x_1+d-x)!}\cdot e^{-(\theta+x)\mu t}(1-e^{-\mu t})^{x_1+d-x}
  \1\{x_1+d-x\ge 0\}
\end{equation}
represents a density for the probability distribution of $(\bar{J}_1,X(\bar{J}_1),\abs{\Delta 
X(\bar{J}_1)})$, w.r.t.\ the product of the Lebesgue measure on $\R^+$ with the counting measure on $S$ and 
with the counting measure on $S\setminus\{0\}$. 
\end{regtheorem}

As a consequence, it is easy to conclude that
\begin{equation}\label{eq-21}
  \prod_{k=1}^{n}f(x_{k-1};\Delta t_k,x_k,d_k)\1\{x_k+d_k-x_{k-1}\ge 0\}
\end{equation}
represents a density for the probability distribution of $(\bar{J}_1,X(\bar{J}_1),\abs{\Delta X(\bar{J}_1)}, 
\ldots,\bar{J}_n,\allowbreak X(\bar{J}_n),\abs{\Delta X(\bar{J}_n)})$, w.r.t.\ an appropriate product measure 
that, at this stage, one can easily envisage.  

\section{Statistical estimation of the unknown probability intensities}\label{sec-6}

As remarked in the previous section, the subject should be discussed on the basis of the likelihood function 
determined by \eqref{eq-14}. But, as already noted elsewhere, there is reason to present the nontrivial 
developments of this idea in a new specific paper. Then, here one contents oneself with a partial solution -- 
that, however, might be useful -- based on the omission of any piece of information about jump times. Whence, 
the question now is about the construction of a statistical model consistent with this new situation. To 
start with, consider \eqref{eq-14} for $n=1$, i.e. 
\[
  \bfP^{(\tau)}(\abs{\Delta X(t_1)}=d_1\mid H_{t_1})=\frac{1}{\sum_{d\ge 1}D^{(\tau)}_1(t_1,d_1)}
  D^{(\tau)}_1(t_1,d_1)
\]
which, according to the new scenario, ought to be independent of $t_1$ to represent the probability that, 
conditionally on the occurrence of a generic negative jump, its magnitude be $d_1$. Hence, that probability 
should coincide with the limit as $t_1\to+\infty$, that is
\begin{multline}\label{eq-22}
  \varphi_\theta(d_1)=\frac{1}{\sum_{n\ge 1}{\theta^n}/{(\theta+1)^n}}\cdot
  \frac{\theta^{d_1}}{(\theta+1)_{d_1}} \\[.5em]
  =\frac{1}{\Phi(1,\theta+1,\theta)-1}\cdot
  \frac{\theta^{d_1}}{(\theta+1)_{d_1}}\qquad\text{($d_1=1,2,\ldots$)}.
\end{multline}
Moreover, the new scenario implies that magnitudes cannot be distinguished from one another, implying that 
they may be seen as exchangeable or, more specifically, conditionally i.i.d.\ according to \eqref{eq-22}. 
Then, the ensuing statistical model $\mathcal{M}$, which depends on $\theta$ only, can be described in the 
following way 
\[
  \mathcal{M}\coloneq \{(D^{\infty},\mathcal{P}(D)^\infty), \varphi^\infty_\theta : \theta>0\}
\]
where: $D\coloneq \Z^+\setminus\{0\}$; $\varphi_\theta(d)$ is the same as \eqref{eq-22} with $d_1=d\in D$. 
This rewriting is done in order to deal with the statistical estimation of $\theta$ as clearly as possible. 

Expectation, variance and other moments of $\varphi_\theta$ are computed in \ref{app-A1} and approximated in 
\ref{app-A4}. The main points can be summarised by saying that the expectation $m(\theta)$ is given by 
\begin{gather*}
  m(\theta) \coloneq \sum_{n\ge 1}n\varphi_\theta(n)=\frac{\theta}{\Phi(1,\theta+1,\theta)-1} \qquad
  \text{($\theta>0$)} \\[.5em]
   \sim \sqrt{\frac{2}{\pi}\theta} \qquad\text{($\theta\to+\infty$)}
\end{gather*}
and the variance by
\begin{gather*}
  \sum_{n\ge 1}n^2\varphi_\theta(n)-m(\theta)^2 =\theta+m(\theta)[1-m(\theta)]\qquad
  \text{($\theta>0$)} \\[.5em]
   \sim \frac{\pi-2}{\pi}\theta  \qquad\text{($\theta\to+\infty$)}.
\end{gather*}

\subsection{Definition of the estimator}\label{sub-6.1}

The estimator in question is a moment estimator as solution to the equation
\[
  \varphi^\infty_\theta(\tilde{d}_1)=\frac{1}{n}\sum_{k=1}^{n}\tilde{d}_k \eqcolon \bar{d}^{(n)}
\]
where $n$ is a fixed positive integer, and $\tilde{d}_1,\tilde{d}_2,\ldots$ are the random coordinate maps of 
$D^\infty$ (i.i.d.\ observations, in the usual statistical language). In view of the expression of 
$m(\theta)$ ($=\varphi^{\infty}_\theta(\tilde{d}_1)$), the above equation becomes 
\begin{equation}\label{eq-23}
  \Phi(1,\theta+1,\theta)=1+\frac{\theta}{\bar{d}^{(n)}}.
\end{equation} 

\begin{regtheorem}\label{th-J}
  Let $\theta_0>0$ denote the ``true'' value of the unknown parameter $\theta$,
  \begin{multline*}
    D_0\coloneq\Bigl\{(d_1,d_2,\ldots)\in D^\infty:\exists n_0=n_0(d_1,d_2,\ldots) \textup{ such that } 
    \textstyle{\sum_{k=1}^{n}}d_k/n>1\\ 
    \textup{ for every } n\ge n_0\Bigr\}.
  \end{multline*}
  Then, $\varphi^\infty_{\theta_0}(D_0)=1$. Moreover, equation~\eqref{eq-23}, to be solved for $\theta>0$, has a unique 
  solution on $D_0$ whenever $n\ge n_0$.
\end{regtheorem}

\begin{proof}[Proof of \textup{(\ref{th-J})}]
  The first assertion is a straightforward consequence of the fact that\linebreak $\varphi^\infty_\theta(\bigcap_{k\ge 1}\{
  \tilde{d}_k=1\})=0$ combined with the observation that $\bar{d}^{(n+k)}>1$ for every $k$, whenever 
  $\bar{d}^{(n)}>1$. To prove the latter, it is enough to highlight a few qualitative features of
  \[
    l(\theta)\coloneq \Phi(1,\theta+1,\theta)=\sum_{n\ge 0}\frac{\Gamma(\theta+1)}{\Gamma(\theta+1+n)}\theta^n \qquad \text{($\theta>0$)}.
  \]
  First of all, from well-known properties of the logarithmic derivative of the gamma function (the psi, or 
  digamma function),
  \begin{align*}
    \frac{d}{d\theta} \frac{\Gamma(\theta+1)}{\Gamma(\theta+1+n)}\theta^n & =\frac{\Gamma(\theta+1)}{\Gamma(\theta+1+n)}\theta^n
    \biggl(\frac{n}{\theta}-\sum_{k=1}^{n}\frac{1}{\theta+k}\biggr) && \text{($n\ge 1$)} \\[.5em]
    & <n\frac{\Gamma(\theta+1)}{\Gamma(\theta+1+n)}\theta^{n-1}<\frac{A^{n-1}}{\Gamma(n)} && \text{($0<\theta<A$)}.
  \end{align*}
  Hence, by the Weierstrass M-test, the series of the derivatives converges uniformly on $(0, A)$ for every 
  $A>0$, which entails
  \[
    l'(\theta)=\sum_{n\ge 1}\frac{\Gamma(\theta+1)}{\Gamma(\theta+1+n)}\theta^{n-1}\sum_{k=1}^{n}\frac{k}{\theta+k} \qquad
    \text{($\theta>0$)}.
  \]
  Then, $l$ is strictly increasing. In a similar way one gets
  \begin{align*}
    l''(\theta) & =\sum_{n\ge 1}\frac{\Gamma(\theta+1)}{\Gamma(\theta+1+n)}\theta^n\Biggl[\biggl(\,\sum_{k=1}^n\frac{k}{\theta(\theta+k)}
    \biggr)^{\!\!2}-\sum_{k=1}^{n}\frac{k^2+2\theta k}{(\theta(\theta+k))^2}\Biggr] \\[.5em]
     & <\sum_{n\ge 1}\frac{\Gamma(\theta+1)}{\Gamma(\theta+1+n)}\theta^n\Biggl[\biggl(\,\sum_{k=1}^n\frac{k}{\theta(\theta+k)}
    \biggr)^{\!\!2}-\sum_{k=1}^{n}\biggl(\frac{k}{\theta(\theta+k)}\biggr)^{\!\!2}\Biggr] \\[.5em]
     & <0
  \end{align*}
  where the last inequality follows from an application to the difference in square brackets of Theorem~19 in 
  \cite{Hardy1952}. To complete the proof, extend the definition of $l$ to $\R^+$, by continuity, to 
  obtain $l(0)=1$, $l'_+(0)=1$, and observe that $y$-intercept and slope of the right-hand-side of 
  \eqref{eq-23} are equal to $1$ and $1/\bar{d}^{(n)}$, respectively, with $\bar{d}^{(n)}>1$ on $D_0$.
\end{proof}
 
The extended equation to $\R^+$ has, therefore, a unique solution $\hat{\theta}_n$, which is equal to $0$ if 
and only if $\bar{d}^{(n)}=1$. According to (\ref{th-J}), the sequence $(\hat{\theta}_n)_{n\ge 1}$ may have 
an initial segment of $0$'s, but, with $\varphi^\infty_{\theta_0}$-probability $1$, there is a random 
positive integer $\tilde{n}_0$ such that $\hat{\theta}_n>0$ for every $n\ge \tilde{n}_0$. Any sequence of 
this type will be called \emph{sequence of moment estimators of $\theta$}. 
 
\subsection{Asymptotic properties of the estimator}\label{sub-6.2}

\begin{figure}[htbp!]
    \centering
    \subfloat[The function $L$]{\includegraphics[width=.45\linewidth]{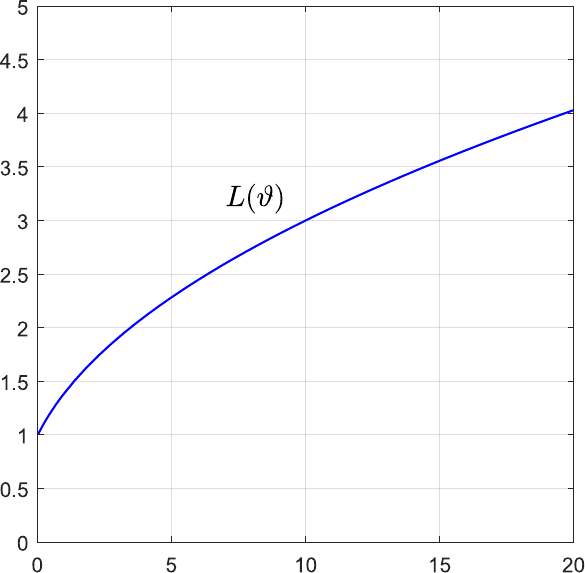}\label{fig-01a}}
    \hfill 
    \subfloat[The function $L'$]{\includegraphics[width=.45\linewidth]{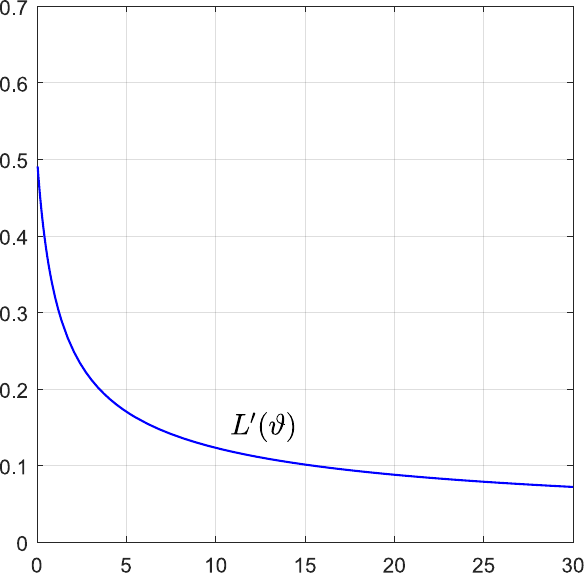}\label{fig-01b}}
    \par\vspace{3em}
    \subfloat[The function $L''$]{\includegraphics[width=.45\linewidth]{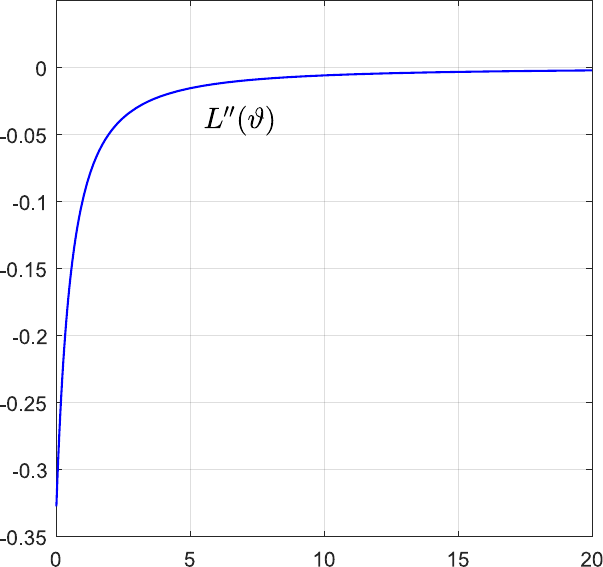}\label{fig-01c}}
    \medskip
    \caption{Graphs of the functions $L$, $L'$ and $L''$.}\label{fig-01}
\end{figure}

In order to establish consistency and asymptotic normality of $(\hat{\theta}_n)_{n\ge 1}$, it is worth 
referring to the expression of $\varphi^\infty_\theta(\tilde{d}_1)$ given in~\ref{app-A1}, i.e.,
\[
  \varphi^\infty_\theta(\tilde{d}_1)=\frac{\theta}{\Phi(1,\theta+1,\theta)-1}=\frac{\theta}{l(\theta)-1} 
  \qquad \text{($\theta>0$)}
\]
which admits the following \emph{continuous} extension to $\R^+$
\[
  L(\theta)=\begin{cases}
              1 & \text{if $\theta=0$} \\[.5em]
              \dfrac{\theta}{l(\theta)-1} & \text{if $\theta>0$}.
            \end{cases}
\]
$L\colon\R^+\to[1,+\infty)$ is \emph{bijective}, as one can ascertain by observing that: 
$L(\theta)\to+\infty$ as $\theta\to+\infty$, from \ref{app-A4}; $L$ is differentiable everywhere, 
$L'_+(0)=1/2=\lim_{\theta\downarrow 0}L'(\theta)$; the derivative $L'(\theta)$ is strictly positive for every 
$\theta>0$. To verify this last fact, suppose there is $\theta_1>0$ for which $L'(\theta_1)=0$. Then, 
\[
  l'(\theta_1)=\frac{l(\theta_1)-1}{\theta_1}=\frac{l(\theta_1)-l(0)}{\theta_1}
\]
and -- in view of the Cavalieri-Lagrange mean-value theorem -- there is $\theta_2$ in $(0,\theta_1)$ such 
that $l'(\theta_1)=l'(\theta_2)$, contradicting the fact that $l''(\theta)<0$ for every $\theta>0$ (see 
computations in the proof of (\ref{th-J})). Further information about $L$ can be found in \ref{app-A2}. For 
the sake of illustration, the graphs of $L$, $L'$ and $L''$ are shown in Figure~\ref{fig-01}. The way is now 
paved for the proof that $(\hat{\theta}_n)_{n\ge 1}$ is strongly consistent. 

\begin{regtheorem}\label{th-K}
  Let $\theta_0>0$ be the true value of the unknown parameter, $(\hat{\theta}_n)$ a sequence of moment 
  estimators of $\theta$ \textup(see equation~\textup{\eqref{eq-23}}\textup), $m$ and $\epsilon$ strictly positive elements 
  of $\Z^+$ and $\R^+$, respectively. Then
  \[
    \varphi^\infty_{\theta_0}\Bigl(\,\bigcap\nolimits_{k\ge m}\{\abs{\hat{\theta}_k-\theta_0}<\epsilon\}\Bigr) \ge
    1-\biggl(\frac{v(\theta_0)}{m(\theta_0,\epsilon)}\biggr)^{\!\!2}\biggl(\frac{1}{m}+
    \sum_{k\ge m+1}\frac{1}{k^2}\biggr)
  \]
  where
  \[
    v^2(\theta)\coloneq\varphi^\infty_\theta((\hat{d}_1-L(\theta))^2)=\theta+L(\theta)(1-L(\theta))
  \]
  and
  \[
    m(\theta,\epsilon)\coloneq L(\theta+\epsilon)-L(\theta)
  \]
  are defined for every $\theta>0$.
\end{regtheorem}

\begin{proof}[Proof of \textup{(\ref{th-K})}]
  Note that the event $\{\abs{\bar{d}^{(n)}-L(\theta)}<m(\theta,\epsilon)\}$ entails $\{\abs{\hat{\theta}_n 
  -\theta_0}<\epsilon\}$, for every $n$, $\theta$, $\epsilon$. Whence the thesis obtains by the Hájek-Rényi 
  maximal inequality (see, e.g., \cite{Petrov1995}, Theorem~2.5) applied to the sequence 
  $(\bar{d}^{(n)})_{n\ge 1}$.
\end{proof}

The above-mentioned characteristics of $L$ also allow application of the \emph{delta method} to obtain

\begin{regtheorem}\label{th-L}
  Let the model $\mathcal{M}$ defined in Subsection~\textup{\ref{sub-6.2}} be in force and $\theta_0$ the ``true'' value of the unknown parameter~$\theta$. 
  Then, $(\sqrt{n}(\hat{\theta}_n-\theta_0))_{n\ge 1}$ converges in law to a Gaussian random 
  quantity with zero expectation and variance $(v(\theta_0)/L'(\theta_0))^2$.
\end{regtheorem}

\begin{proof}[Proof of \textup{(\ref{th-L})}]
  By conjunction of the delta method with the central limit theorem for 
  $\bigl((\bar{d}^{(n)}-L(\theta_0))\sqrt{n}\,\bigr)_{n\ge 1}$.
\end{proof}

\subsection{Brief mention to the estimation of the rate of negative jumps}\label{sub-6.3} 

The conditions under which estimation has been tackled so far imply that $\mu$ -- the probability intensity 
of a single negative jump of any magnitude -- be ignored. All the same, the estimation of $\mu$ may benefit 
from the fact that it is directly connected with the sole facts one can actually observe. In point of fact, 
one could consider, as an estimator of $\mu$, the empirical average number $\hat{\mu}$ of negative jumps per 
unit of time, say $u$ (e.g., the day, the week, \ldots). Then, like in \cite{deFinetti1957}, by adopting 
$u/\hat{\mu}$ as new operating unit of time, one can assume $\mu=1$ and $\theta=\lambda$. 

\subsection{Final remarks from a Bayesian perspective}\label{sub-6.4}

It should be noted that both the coefficient $(v/m)^2$ in (\ref{th-K}) and the variance of the limiting law 
in (\ref{th-L}) depend on $\theta_0$, that is a hypothetical value of an unknown parameter. This event, which 
generally affects similar statistical statements of a frequentistic nature, compromises value and practical 
usefulness of them. It stems from the prejudice against the inclusion in the theory of statistics of the 
initial opinions -- in the form of subjective probabilities -- of the ``user'' of the theory. But ``\ldots\ 
no problem can be correctly stated in statistics without an evaluation of the initial probabilities\ldots\ 
[even though]\ldots\ probabilities must often be based on vague, uncertain, and fragmentary information'' 
(\cite{deFinetti1972}, Chapter~8). A situation which could arise also in connection with the problem of 
estimating $\theta\coloneq \lambda/\mu$, because of the supposed little reliable information about $\lambda$. 
In any case, in order to take due account of initial probabilities, one ought to modify and complete the 
model $\mathcal{M}$ according to the usual Bayes-Laplace paradigm. Thus, the space of the ``observables'' is 
replaced with the product of the parameter space and the space of observables, viz, $(\R^+\times 
D^\infty,\mathcal{B}(\R^+)\otimes(\mathcal{P}(D))^\infty)$. Its first coordinate, say $\tilde{\theta}$, 
represents the random parameter, while the subsequent ones, denoted by $\tilde{d}_1,\tilde{d}_2, \ldots$\,, 
have the same meaning as in $\mathcal{M}$.  Denoting by $q$ the initial distribution of $\tilde{\theta}$ -- 
i.e., a probability measure on $(\R^+,\mathcal{B}(\R^+))$ -- the assumption that the elements of 
$\mathcal{M}$ form a conditional distribution for $(\tilde{d}_n)_{n\ge 1}$, given $\tilde{\theta}$, is 
tantamount to endowing the above product space with the probability measure $Q$  presentable as 
\[
  Q(d\theta\, d\mathbf{x})=q(d\theta)\varphi^\infty_\theta(d\mathbf{x}) \qquad \text{($\theta\ge 0$, 
  $\mathbf{x}\in D^\infty$)}
\]
with the proviso that $\varphi_0\coloneq\delta_1$. It is important to note that $\tilde{d}_1,\tilde{d}_2, 
\ldots$ are \emph{exchangeable} w.r.t.\ $Q$ and their (joint) distribution is $\int_{\R^+}\varphi^\infty_0\, 
dq$.

Based on these premises, one can now try to find a sort of Bayesian justification for the moment estimator 
defined in Subsection~\ref{sub-6.1}, and, then, discuss Bayesian counterparts of (\ref{th-K}), (\ref{th-L}), 
which, of course, are free from the defect mentioned at the beginning of the present one. As far as the first 
point is concerned, let $\hat{s}_n=\hat{s}_n(\tilde{d}_1,\ldots,\tilde{d}_n)$ denote any potential estimator 
of $\tilde{\theta}$, and suppose that the loss suffered by the ``user'' for deciding $\hat{s}_n$ is 
$\abs{L(\hat{s}_n)-\bar{d}^{(n)}}$. This loss function, although rather odd, is not totally senseless. Since 
the corresponding statistical risk $Q(\abs{L(\hat{s}_n)-\bar{d}^{(n)}})$ vanishes when 
$\hat{s}_n=\hat{\theta}_n$, then the moment estimator $\hat{\theta}_n$ can be viewed as a Bayesian estimator 
of $\tilde{\theta}$. In the new Bayesian context the study of the consistency of $\hat{\theta}_n$ leads to 
consider $Q(\bigcap_{k\ge m}\{ \abs{\hat{\theta}_k-\tilde{\theta}}<\epsilon\})$ which, by virtue of 
(\ref{th-K}), obeys 
\begin{equation}\label{eq-24}
  Q\Big(\bigcap\nolimits_{k\ge m}\{\abs{\hat{\theta}_k-\tilde{\theta}}<\epsilon\}\Big) \ge
  1-\biggl(\frac{1}{m}+\sum_{k\ge m+1}\frac{1}{k^2}\biggr)
  \int_{\R^+}\biggl(\frac{v(x)}{L(x+\epsilon)-L(x)}\biggr)^{\!\!2}\,q(dx)
\end{equation}
\emph{provided that the integral is finite}. Supporters of non-Bayesian views of statistics might object that 
the Bayesian answer provided by \eqref{eq-24}, to the original critical issue about the indeterminacy of the 
bound in (\ref{th-K}), is simply illusory. Indeed, some say that the evaluation of $q$  is impractical, 
others say that prior distributions do not exist. But, in the present case one should note that the new bound 
boils down to the expectation of 
\[
  W(\tilde{\theta})=\biggl(\frac{v(\tilde{\theta})}{L(\tilde{\theta}+\epsilon)-L(\tilde{\theta})}\biggr)^{\!\!2}
  =\frac{\tilde{\theta}+L(\tilde{\theta})(1-L(\tilde{\theta}))}{(\epsilon\cdot L'(\tilde{c}))^2}
\]
where $\tilde{c}$ is some suitable point of the interval $(\tilde{\theta},\tilde{\theta}+\epsilon)$, and, to 
evaluate it, complete knowledge of $q$ is not required at all. In \ref{app-A3} and \ref{app-A4} there are 
computations of the elements of the above expression for small values and, more importantly, for large values 
of $\tilde{\theta}$, respectively. In fact, for the latter, in \ref{app-A4} it is proven that 
\[
  W(\tilde{\theta})\sim \frac{2(\pi-2)}{\epsilon^2}\bigl(\tilde{\theta}^2+\epsilon\cdot
  c(\epsilon,\tilde{\theta})\tilde{\theta}\bigr) 
\]
where $c(\epsilon,\tilde{\theta})$ is a random number taking value in $(0,1)$. Then, with a view to practical 
uses of \eqref{eq-24}, it suffices that the ``user'' forms an opinion on, and accordingly assesses, the first 
two moments of $q$ directly. From a subjectivistic viewpoint, the sole requirement such an opinion has to 
satisfy is that of being coherent w.r.t.\ the rest of her/his initial opinions. 

Reconverting (\ref{th-L}) to a Bayesian statement seems to be of little interest and even contradictory. 
Indeed, the assumption that $q$ has finite second moment is sufficient in order that, by exchangeability, 
$(\bar{d}^{(n)})_{n\ge 1}$ converges to a random number $\tilde{d}$, a.s.-$Q$. To determine the corresponding 
(limiting) distribution of $\hat{\theta}_n$, consider any open set $G$ of $\R^+$, meant as topological 
subspace of $\R$ endowed with the usual metric. Then, 
\begin{align*}
  \liminf_{n\to+\infty} Q\{\bar{d}^{(n)}\in G\} & =\liminf_{n\to+\infty}\int_{\R^+}\varphi_\theta^\infty\{ 
  \bar{d}^{(n)}\in G\}\,q(d\theta) \\[.5em]
    & \ge \int_{\R^+}\liminf_{n\to+\infty}\varphi_\theta^\infty\{\bar{d}^{(n)}\in G\}\,q(d\theta) && 
    \text{(by Fatou's lemma)} \\[.5em]
    & \ge \int_{\R^+}\delta_{L(\theta)}(G)\,q(d\theta)=q(L^{-1}(G))
\end{align*}
where the second inequality follows from the conjunction of the Kolmogorov strong law of large numbers with 
the necessary part of Theorem~2.1(iv) -- portmanteau theorem -- in \cite{Billingsley1999}. Whence, the 
sufficient part of that very same theorem implies that $\tilde{d}$ must be distributed according to $q\circ 
L^{-1}$. Then, by the (continuous mapping)~Theorem~2.7 therein, the distribution of  $L^{-1}(\tilde{d}\,)$, 
that is the initial $q$, must be the limiting distribution of the estimator  
$\hat{\theta}_n=L^{-1}(\bar{d}^{(n)})$, as $n\to+\infty$. 

\appendix

\renewcommand{\thesubsection}{A\arabic{subsection}} 

\section*{Appendix}\label{sec-appendix}
\addcontentsline{toc}{section}{Appendix}

\counterwithin*{equation}{subsection}
\renewcommand{\theequation}{\thesubsection.\arabic{equation}}

\allowdisplaybreaks[1] 

Gathered here for the sake of completeness are certain calculations and remarks of a numerical nature 
pertaining to results shown or used in previous sections, specifically the sixth one. They consist in rather 
simple adaptations of well known facts concerning the \emph{incomplete gamma functions} and its asymptotic 
expansion. For easy and exhaustive reference, the reader is referred to Chapters~8 and 13 of 
\cite{Olver2010}, mentioned as \HMF\ throughout this Appendix. It should be noted that the Kummer function 
$\Phi$ is denoted therein by $M$. 

\subsection{Expectation, variance and other moments of the distribution $\varphi_\theta$}\label{app-A1}

If $\theta>0$, then 
\begin{align*}
  \varphi_\theta(\tilde{d}_1) & =\frac{1}{\Phi(1,\theta+1,\theta)-1} 
                                \sum_{d\ge 1}\frac{\theta^d}{(\theta+1)_d} \\[.5em]
                              & =\frac{\theta}{\Phi(1,\theta+1,1)-1} \frac{d}{dz}
                                \,\Phi(1,\theta+1,z)\Big\lvert_{z=\theta} \\[.5em]
                              & =\frac{\theta}{\Phi(1,\theta+1,1)-1}\bigl[\theta(\theta+1)+(\theta+1)(z-\theta)
                                 \Phi(1,\theta+1,z)\bigr]\Big\lvert_{z=\theta}
\end{align*}
where the last equality follows from \HMF~13.3.15 and recurrence relations 13.3. Then
\[
  \varphi_\theta(\tilde{d}_1)=\frac{\theta}{\Phi(1,\theta+1,\theta)-1} \qquad \text{($\theta>0$)}
\]
expression already used in Subsection~{\ref{sub-6.2} and denoted by $L(\theta)$ therein.

In the same vein, given any integer $m\ge 2$, by resorting to the horizontal generating function of the 
Stirling numbers of the second kind $S(n,k)$ -- see \HMF~26.8.10 -- one gets
\[
  \varphi_\theta(\tilde{d}^{\,m}_1)=\frac{1}{\Phi(1,\theta+1,\theta)-1}\sum_{k=1}^m \theta^k S(m,k)
  \Phi(1+k,\theta+1+k,\theta).
\]
In particular, for $m=2$, the aforesaid recurrence relations can be applied to the above expression to obtain
\[
  \varphi_\theta(\tilde{d}^{\,2}_1)=\theta+L(\theta) \qquad \text{($\theta>0$)}
\]
and then
\[
  v(\theta)^2\coloneq \varphi_\theta(\tilde{d}^{\,2}_1)-\varphi_\theta(\tilde{d}_1)^2 
  =\theta + L(\theta)[1-L(\theta)].
\]

\subsection{Extension of $L$ and its connection with the incomplete gamma 
function}\label{app-A2}

For each natural $n$, $z\mapsto z^n/(z+1)_{n+1}$ is meromorphic with poles at $-1,\ldots,-(n+1)$, i.e., 
analytic on the region of regularity $\C\setminus\{-1,\ldots,-(n+1)\}$. With the help of the representation
\[
  \frac{1}{(z+1)_{n+1}}=\frac{1}{(n+1)!}\exp\biggl\{-\sum_{j=1}^{n+1}\Log\biggl(1+\frac{z}{j}\biggr)\biggr\}
\]
(where $\Log$ denotes the principal value) and the Weierstrass test theorem, it is easy to see that 
$\sum_{n\ge 0}z^n/(z+1)_{n+1}$ is uniformly convergent on $\mathcal{R}\coloneq\C\setminus\{-1,-2,\ldots\}$. 
Then, by the Weierstrass theorem about series of analytic functions, one establishes that 
\[
  D(z)\coloneq \frac{l(z)-1}{z}\coloneq \sum_{n\ge 0}\frac{z^n}{(z+1)_{n+1}}
\]
is analytic on $\mathcal{R}$. As a consequence, if $\mathcal{Z}$ denotes the set of the zeros of $D$, then 
$L(z)\coloneq 1/D(z)$ is regular on $\mathcal{R}_0\coloneq\mathcal{R}\setminus\mathcal{Z}$. When $z=\theta\ge 
0$, $L$ is the same as the namesake function already defined in Subsection~\ref{sub-6.2}. One now derives the 
expansion of $D$ into powers of $z$, about $z=0$, from an identity by \cite{Cifarelli1979a,Cifarelli1979b} -- 
see Theorem~1 in \cite{Lijoi2004} -- which entails 
\begin{equation}\label{eq-a2-1}
\begin{split}
  \frac{1}{(z+1)_{n+1}} & =\frac{1}{(n+1)!}\exp\biggl\{-\sum_{j=1}^{n+1}\Log\biggl(1+\frac{z}{j}\biggr)\biggr\} 
                            \\[.5em]
                        & =\frac{1}{(n+1)!}\int_{\mathcal{C}^{(n)}}(1+z\mathcal{H}^{-1}_n(p))^{-(n+1)}
                        \,\mathcal{D}_n(dp)
\end{split}
\end{equation}
where: $\mathcal{C}^{(n)}$ is the set of all probabilities concentrated on $\{1,\ldots,n+1\}$, 
$\mathcal{H}^{-1}_n(p)\coloneq p_1+p_2/2+\ldots+ p_{n+1}/(n+1)$, i.e., the reciprocal of the harmonic mean 
$\mathcal{H}_n(p)$  of the probability $p$ on $\{1,\ldots,n+1\}$, and $\mathcal{D}_n$ the Ferguson-Dirichlet 
distribution on $\mathcal{C}^{(n)}$  with parameter given by the counting measure on $\{1,\ldots,n+1\}$. If 
$\abs{z}<1$, then also $\abs{z}\mathcal{H}^{-1}_n(p)<1$ and, by the binomial series theorem, one obtains 
\[
  \frac{1}{(z+1)_{n+1}}=\frac{1}{(n+1)!}\int_{\mathcal{C}^{(n)}}\sum_{k\ge 0}(-1)^k\frac{(n+1)_k}{k!}
  (z\mathcal{H}^{-1}_n(p))^k\,\mathcal{D}_n(dp)
\]
where
\begin{multline*}
  \int_{\mathcal{C}^{(n)}}\sum_{k\ge 0}\frac{(n+1)_k}{k!}\abs{z}^k(\mathcal{H}^{-1}_n(p))^k\,\mathcal{D}_n(dp) \\[.5em]
  \leq \int_{\mathcal{C}^{(n)}}\sum_{k\ge 0} \frac{(n+1)_k}{k!}\abs{z}^k\,\mathcal{D}_n(dp) 
  \eqcolon {}_2F_1(a,n+1;a;\abs{z})
\end{multline*}
$a$ is any positive, and, as usual, ${}_2F_1$ denotes the \emph{Gaussian} (or \emph{ordinary}) 
\emph{hypergeometric series}. Thus, one gets 
\[
  \frac{1}{(z+1)_{n+1}}=\frac{1}{(n+1)!}\sum_{k\ge 0}(-1)^k \frac{(n+1)_k}{k!} \mu^{(n)}_k z^k \qquad 
  \text{($\abs{z}<1$)}
\]
where $\mu^{(n)}_k\coloneq \int_{\mathcal{C}^{(n)}}(\mathcal{H}^{-1}_n(p))^k\,\mathcal{D}_n(dp)\in(0,1]$, 
($k=0,1,\ldots$), and
\begin{equation}\label{eq-a2-2}
  D(z)=\frac{l(z)-1}{z}=\sum_{n\ge 0} \frac{z^n}{(n+1)!}\sum_{k\ge 0}(-1)^k\frac{(n+1)_k}{k!}\mu^{(n)}_k z^k.
\end{equation} 
The series of the absolute values may be majorised by
\[
  \sum_{n\ge 0}\frac{\abs{z}^n}{(n+1)!}\,{}_2F_1(a,n+1;a;\abs{z})
\]
which, since ${}_2F_1(a,n+1;a;\abs{z})\sim e^{(n+1)\abs{z}}$ as $n\to+\infty$ (see 2.3.2(13) in 
\cite{Erdely1953}), converges. This allows one to interchange the order of summation in \eqref{eq-a2-2} to 
obtain
\[
  \frac{l(z)-1}{z}=\sum_{\nu\ge 0}c_\nu z^\nu \qquad \text{($\abs{z}<1$)}
\]
where
\begin{equation}\label{eq-a2-3}
  c_\nu=\sum_{k=0}^{\nu}(-1)^k\,\frac{(\nu-k+1)_k}{(\nu-k+1)!} \, \frac{\mu_k^{(\nu-k)}}{k!}
\end{equation}
for $\nu=0,1,\ldots$\;. To compute $\mu^{(n)}_k$ one may equate the $k$-th derivatives, at $z=0$, of both 
sides of \eqref{eq-a2-1} (for the left side one might need the Faà di Bruno formula concerning the successive 
derivatives of composite functions). Since $D(0)=1$, there is $r\in (0,1)$ for which the following expansion 
of $L$ can be derived by means of the ordinary rules, i.e. 
\[
  L(z)=\sum_{n\ge 0} \beta_n z^n \qquad \abs{z}<r
\]
with
\begin{align*}
  \beta_0 & =1/c_0=1 \\
          & \vdots \\
  \beta_n & =-\sum_{\nu=1}^{n} c_\nu\beta_{n-\nu} \qquad \text{($n=1,2,\ldots$)}
\end{align*}
$c_0,c_1,\ldots$ being the same as in \eqref{eq-a2-3}.

For the sake of illustration, one determines $\beta_1$ and $\beta_2$. 
\begin{align*}
  \beta_1=-c_1\beta_0=-c_1 & =-\sum_{k=0}^{1}(-1)^k\, \frac{(2-k)_k}{(2-k)!}\,\frac{\mu_k^{(1-k)}}{k!} \\[.5em]
    & =-\biggl(\frac{1}{2}-1\biggr)=\frac{1}{2}
\end{align*}
entailing $L'(0)=\beta_1=1/2$;
\begin{align*}
  \beta_2=-\sum_{\nu=1}^{2}c_\nu\beta_{2-\nu} & =-\biggl(-\frac{1}{4}+c_2\biggr) \\[.5em]
                 & =\frac{1}{4}-\sum_{k=0}^{2}(-1)^k\,\frac{(3-k)_k}{(3-k)!}\,\frac{\mu^{(2-k)}_k}{k!} \\[.5em]
                 & =\mu^{(1)}_1-\frac{11}{12} \\[.5em]
                 & =\frac{1}{2}\biggl(1+\frac{1}{2}\biggr)-\frac{11}{12} \\[.5em]
                 & =-\frac{1}{6}
\end{align*}
entailing $L''(0)=2\beta_2=-1/3$. Cf.\ Figures~\ref{fig-01b} and \ref{fig-01c}.

The function $D$ can also be seen as analytic continuation of a remarkable analytic special function defined 
on $\{z\in\C:\RE z>0\}$. In fact, from a classical result (for this and next developments, cf.~\HMF, 
Chapter~8), 
\[
  \Phi(1,1+a,z)=a z^{-a} e^z \gamma(a,z) \qquad \text{($\RE a>0$)}
\]
where $z^{-a}$ stands for the principal value of the power, and
\[
  \gamma(a,z)\coloneq \int_0^z x^{a-1}e^{-x}\, dx
\]
is just the incomplete gamma function with the restriction that the integration path does not cross the 
negative real axes. Then, 
\begin{align*}
  D(z) & =e^z\cdot z^{-(z+1)}\gamma(z+1,z) \qquad \text{($\RE z>0$)} \\[.5em]
       & (=\Gamma(z+1)e^z z^{-(z+1)}P(z+1,z))
\end{align*}
and, by virtue of a well-known recurrence relation,
\[
  z\cdot D(z)=z\Gamma(z)\biggl(\frac{e}{z}\biggr)^{\!\!z} P(z,z)-1 \qquad \text{($\RE z>0$)}.
\]

These remarks will come in handy in \ref{app-A3} and \ref{app-A4}, in which the problem of finding exact and 
asymptotic evaluations of $L(\theta)$, and related functions, is dealt with. 

\subsection{Evaluation of $L(\theta)$ for small and moderate values of $\theta$}\label{app-A3}

With a view to computing $L(\theta)$ it suffices to pay attention to $D(\theta)$ which, according to the 
final statements of \ref{app-A2}, satisfies 
\begin{align*}
  \theta D(\theta) & =\theta \Gamma(\theta)\biggl(\frac{e}{\theta}\biggr)^{\!\!\theta} P(\theta,\theta)-1
  \qquad \text{($\theta>0$)} \\[.5em]
   & =\theta\Gamma(\theta) e^{\theta}\gamma^*(\theta,\theta)-1
\end{align*} 
where $\gamma^*$ is the well-known Tricomi variant of $\gamma$. If $\theta$ is not large, then one may 
compute $D$ by resorting to certain numerical series deduced from 8.7.1 in \HMF, to obtain 
\begin{align*}
  D(\theta) & =\Gamma(\theta)\sum_{k\ge 1}\frac{\theta^k}{\Gamma(\theta+k+1)} \\[.5em]
     & =e^{\theta}\sum_{k\ge 1}\frac{(-1)^k}{k!(\theta+k)} \qquad \text{($\theta>0$)}.
\end{align*}

There are also various tables of values of $\gamma,P,\ldots$\;. It is worth mentioning the classical work of 
\cite{Pearson1922} in which the values of 
\[
  I(u,p)\coloneq \frac{1}{\Gamma(p+1)}\int_{0}^{u\sqrt{p+1}} e^{-x}x^p\, dx
\]
are tabled for $p>-1$ and $u>0$, $I$ being obviously connected with $P$ through the relation
\[
  P(\theta,x)=I\biggl(\frac{x}{\sqrt{\theta}},\theta-1\biggr) \qquad \text{($\theta,x>0$)}.
\]

Of course, there is also a great deal of software for computing the above functions.

\subsection{Asymptotic approximation of $L$ (large $\theta$)}\label{app-A4}

The representations of $D$ displayed at the end of \ref{app-A2} are likewise suitable for the study of the 
asymptotic behaviour of $D$ and related functions, as $z\to\infty$. Before dealing with the rather particular 
cases of interest for the present paper, it is worth recalling that the assertion ``the formal series 
$\sum_{s=0}^\infty a_s z^{-s}$ is an \emph{asymptotic} (or \emph{Poincaré}) expansion of the function $f(z)$ 
-- in symbols, $f(z)\sim\sum_{s\ge 0}a_s z^{-s}$ -- as $z\to\infty$ in the unbounded set 
$\mathcal{C}\subset\C$'', means that 
\begin{equation}\label{eq-a4-1}
  z^n\biggl\{f(z)-\sum_{s=0}^{n-1}\frac{a_s}{z^s}\biggr\}\to a_n
\end{equation}
as $z\to\infty$ in $\mathcal{C}$, uniformly w.r.t.\ $\arg z$, for each natural~$n$.

In point of fact, it follows from \HMF~8.12.18--8.12.20 and \cite{Paris2016} that
\begin{equation}\label{eq-a4-2}
  \sqrt{\theta}D(\theta)\sim\sum_{s\ge 0}f_s \theta^{-s/2} \qquad \text{($\theta\to+\infty$)}
\end{equation}
where
\[
  f_0=\sqrt{\frac{\pi}{2}},\quad f_1=-\frac{2}{3},\quad f_s=\sqrt{\frac{\pi}{2}}A_s+B_s 
  \quad\text{for $s=2,3,\ldots$}
\]
and $A_s=A_s(0)$, $B_s=B_s(0)$ are the same as in \cite{Paris2016}, Table~2 with $x=0$: e.g., $A_0=1$, 
$A_2=1/12$, $A_4=1/288$, \ldots\;, $A_1=A_3=A_5=\ldots=0=B_0=B_2=B_4=\ldots$\;, $B_1=1/3$, $B_3=4/135$, 
$B_5=-8/2835$, \ldots\;; hence $f_0\approx 1.25331$, $f_1\approx -0.66667$, $f_2\approx 0.10444$, $f_3\approx 
0.02963$, $f_4\approx 0.00435$, $f_5\approx -0.00282$. 

Consequently, one has 
\[
  \frac{1}{\sqrt{\theta}D(\theta)}\sim \sum_{s\ge 0}\frac{k_s}{\theta^{s/2}} \qquad \text{($\theta\to+\infty$)}
\]
where $k_0=1/f_0=\sqrt{2/\pi}\approx 0.79788$, $k_1=-f_1/f_0^2=4/(3\pi)\approx 0.42441$, $k_2=(f_1^2-f_0 
f_2)/f_0^3=(2/\pi)^{3/2}(4/9-\pi/24)\approx 0.15926$, 
$k_3=(-f_1^3+2f_0f_1f_2-f_0^2f_3)/f_0^4=-(2/\pi)^2\cdot(8/27+\pi(1/18+2/135))\approx -0.20968$, \ldots 
\[
  k_s=-\sqrt{\frac{2}{\pi}}\sum_{i=1}^{s} f_i k_{s-i} \qquad \text{($s=1,2,\ldots$)}
\]
and, more importantly,
\begin{align*}
  L(\theta) & \sim \sqrt{\theta}\sum_{s\ge 0}\frac{k_s}{\theta^{s/2}} \qquad \text{($\theta\to+\infty$)} \\[.5em]
   & \approx \sqrt{\theta}\,\biggl\{0.79788+\frac{0.42441}{\sqrt{\theta}}+\frac{0.15926}{\theta} - 
   \frac{0.20968}{\theta^{3/2}}+\ldots\biggr\}.
\end{align*}

The above asymptotic expansion \eqref{eq-a4-2} is valid for $\abs{z}\to+\infty$ in $\abs{\arg z}<\pi/2$. This 
implies that there are $\rho>0$ and $\eta\in (0,z/2)$ such that $zD(z)$ has no zeros in the annular sector 
$S_{\eta,\rho}\coloneq\{z\in\C:\abs{z}\ge\rho,\abs{\arg z}<\pi/2-\eta\}$. Moreover, since one knows, from 
\ref{app-A2}, that $D$ is analytic on $\C\setminus\{0,-1,-2,\ldots\}$, then $L(z)=\dfrac{z}{zD(z)}$ is 
analytic on $S_{\eta,\rho}$. This fact allows differentiation of the asymptotic expansion of $L$ any number 
of times in the sector $S_{\eta,\rho}$. See \cite{Henrici1977}, Theorem~11.3f. In particular, one obtains 
\[
  L'(\theta)\sim\frac{1}{2}\frac{k_0}{\theta^{1/2}}-\sum_{s\ge 3} \frac{s-2}{2}\frac{k_{s-1}}{\theta^{s/2}}
  \qquad \text{($\theta\to+\infty$)}
\]
and
\[
  L''(\theta)\sim-\frac{1}{4}\frac{k_0}{\theta^{3/2}}+\sum_{s\ge 5} 
  \frac{(s-2)(s-4)}{4}\frac{k_{s-3}}{\theta^{s/2}} \qquad \text{($\theta\to+\infty$)}.
\]
These approximations are visualised in Figure~\ref{fig-02}, juxtaposed to their exact counterparts. 

\begin{figure}[!htbp]
    \centering
    \subfloat[The function $L$]{\includegraphics[width=.45\linewidth]{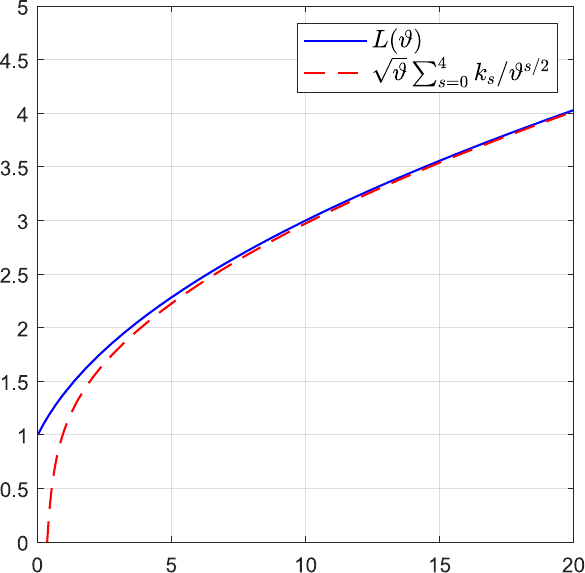}}
    \hfill 
    \subfloat[The function $L'$]{\includegraphics[width=.45\linewidth]{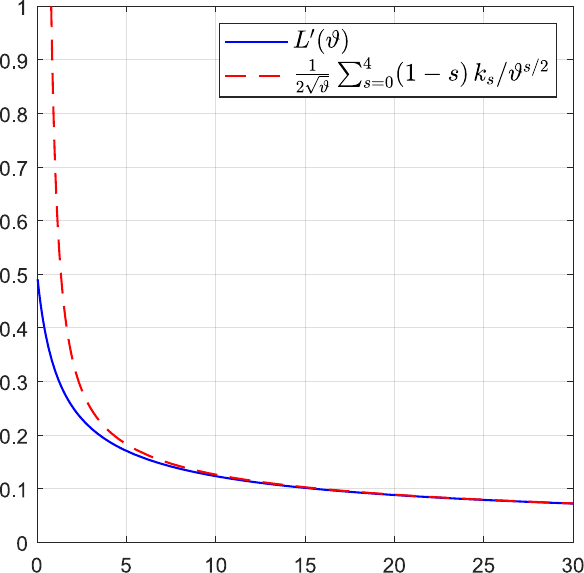}}
    \par\vspace{3em}
    \subfloat[The function $L''$]{\includegraphics[width=.45\linewidth]{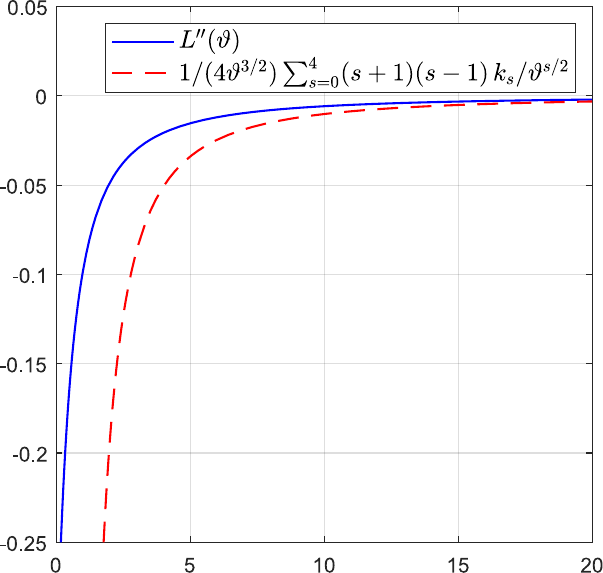}}
    \bigskip
    \caption{Graphs of $L$, $L'$, $L''$ and their asymptotic expansions for 
    $\theta\to+\infty$, truncated at $n=5$, according to~\eqref{eq-a4-1}.}\label{fig-02}
\end{figure}

As an application of the above results, here is an outline of what one should do to obtain an asymptotic 
approximation, as $\theta\to+\infty$, of the quantity $W$ defined in Subsection~\ref{sub-6.4}. The numerator 
reads
\[
  v^2(\theta)=\theta+L(\theta)[1-L(\theta)]\sim \frac{\pi-2}{\pi}\theta+\alpha_1\sqrt{\theta}+
  \alpha_2+\sum_{s\ge 1} \frac{\alpha_{s+2}}{\theta^{s/2}} \qquad \text{($\theta\to+\infty$)}
\]
where $\alpha_n=k_{n-1}-\sum_{j=0}^{n}k_j k_{n-j}$ for $n=1,2,\ldots$\;, and hence: 
$\alpha_1=\sqrt{2/\pi}\,\bigl(1-8/(3\pi)\bigr)$, $\alpha_2=(15\pi-16)/(9\pi^2)$, \ldots\;.

As to the denominator of $W$, it equals $\epsilon^2\cdot(L'(\bar{c}))^2$, where 
$\bar{c}\in(\theta,\theta+\epsilon)$ and
\[
  L'(x)^2\sim\sum_{n\ge 2} \frac{b_n}{x^{n/2}} \qquad \text{($x\to+\infty$)}
\]
with $b_n\coloneq\sum_{j=1}^{n-1}l_j l_{n-j}$ for $n=2,3,\ldots$\;; $l_j\coloneq (1-j/2)k_{j-1}$ for 
$j=1,2,\ldots$\;. Then $b_2=1/(2\pi)$, $b_3=0$, $b_4=-(32-3\pi)/(6\pi)^2$, \ldots\;.

\paragraph*{Acknowledgements.} \addcontentsline{toc}{section}{Acknowledgements} 

My special thanks go to Persi Diaconis for his comments. Warm thanks and appreciation go to Riccardo Dossena 
for technical advice and typing, expert and accurate. 

\nocite{*}

\printbibliography[heading=bibintoc]

\end{document}